\documentclass{article}
\usepackage[utf8]{inputenc}
\usepackage{fullpage}
\usepackage{times}
\usepackage{ifthen}   
\usepackage{graphicx,graphics}
\usepackage{setspace}
\usepackage{amssymb,amsmath,amsthm}
\usepackage{dsfont}
\usepackage{color}
\usepackage{wrapfig}
\usepackage{multirow}
\usepackage{multicol}
\usepackage{pdfpages}
\usepackage{longtable} 	
\usepackage{booktabs}
\usepackage{datetime}
\usepackage{xspace}
\usepackage{paralist} 
\usepackage{enumitem} 
\usepackage[it,small,belowskip=-15pt]{caption}
\usepackage[top=1in,bottom=1in,left=0.95in,right=0.95in]{geometry} 
\usepackage[T1]{fontenc}
\usepackage{doi}

\usepackage{epstopdf}
\usepackage[caption=false]{subfig}
\usepackage[numbers,sort&compress]{natbib}
\renewcommand\bibfont{\fontsize{10}{12}\selectfont}

\theoremstyle{plain}
\newtheorem{theorem}{Theorem}[section]
\newtheorem{lemma}[theorem]{Lemma}
\newtheorem{corollary}[theorem]{Corollary}
\newtheorem{proposition}[theorem]{Proposition}
\theoremstyle{definition}
\newtheorem{definition}[theorem]{Definition}
\newtheorem{example}[theorem]{Example}
\usepackage{hyperref}
\theoremstyle{remark}
\newtheorem{remark}{Remark}
\newtheorem{assumption}{Assumption}

\usepackage{dsfont}
\usepackage{algorithm,algcompatible}

\makeatletter
\let\save@mathaccent\mathaccent
\newcommand*\if@single[3]{%
  \setbox0\hbox{${\mathaccent"0362{#1}}^H$}%
  \setbox2\hbox{${\mathaccent"0362{\kern0pt#1}}^H$}%
  \ifdim\ht0=\ht2 #3\else #2\fi
  }
\newcommand*\rel@kern[1]{\kern#1\dimexpr\macc@kerna}
\newcommand*\widebar[1]{\@ifnextchar^{{\wide@bar{#1}{0}}}{\wide@bar{#1}{1}}}
\newcommand*\wide@bar[2]{\if@single{#1}{\wide@bar@{#1}{#2}{1}}{\wide@bar@{#1}{#2}{2}}}
\newcommand*\wide@bar@[3]{%
  \begingroup
  \def\mathaccent##1##2{%
    \let\mathaccent\save@mathaccent
    \if#32 \let\macc@nucleus\first@char \fi
    \setbox\z@\hbox{$\macc@style{\macc@nucleus}_{}$}%
    \setbox\tw@\hbox{$\macc@style{\macc@nucleus}{}_{}$}%
    \dimen@\wd\tw@
    \advance\dimen@-\wd\z@
    \divide\dimen@ 3
    \@tempdima\wd\tw@
    \advance\@tempdima-\scriptspace
    \divide\@tempdima 10
    \advance\dimen@-\@tempdima
    \ifdim\dimen@>\z@ \dimen@0pt\fi
    \rel@kern{0.6}\kern-\dimen@
    \if#31
      \overline{\rel@kern{-0.6}\kern\dimen@\macc@nucleus\rel@kern{0.4}\kern\dimen@}%
      \advance\dimen@0.4\dimexpr\macc@kerna
      \let\final@kern#2%
      \ifdim\dimen@<\z@ \let\final@kern1\fi
      \if\final@kern1 \kern-\dimen@\fi
    \else
      \overline{\rel@kern{-0.6}\kern\dimen@#1}%
    \fi
  }%
  \macc@depth\@ne
  \let\math@bgroup\@empty \let\math@egroup\macc@set@skewchar
  \mathsurround\z@ \frozen@everymath{\mathgroup\macc@group\relax}%
  \macc@set@skewchar\relax
  \let\mathaccentV\macc@nested@a
  \if#31
    \macc@nested@a\relax111{#1}%
  \else
    \def\gobble@till@marker##1\endmarker{}%
    \futurelet\first@char\gobble@till@marker#1\endmarker
    \ifcat\noexpand\first@char A\else
      \def\first@char{}%
    \fi
    \macc@nested@a\relax111{\first@char}%
  \fi
  \endgroup
}
\makeatother

\DeclareMathOperator*{\argmin}{\arg\!\min}

\newcommand{\real}{{\rm{I\hspace{-.75mm}R}}}

\renewcommand*{~}{\relax\ifmmode\sim\else\nobreakspace{}\fi}

\newcommand{\Var}{\mbox{Var}}

\newcommand{\Dbar}{\bar{D}}
\newcommand{\Ebar}{\bar{E}}
\newcommand{\Fbar}{\bar{F}}

\newcommand{\fhat}{\hat{f}}

\newcommand{\Mhat}{\hat{M}}

\newcommand{\rhohat}{\hat{\rho}}
\newcommand{\sigmahat}{\hat{\sigma}}

\newcommand{\BFe}{\bm{e}}

\newcommand{\BFs}{\bm{s}}

\newcommand{\BFv}{\bm{v}}

\newcommand{\BFx}{\bm{x}}
\newcommand{\BFy}{\bm{y}}

\newcommand{\BFG}{\bm{G}}

\newcommand{\BFX}{\bm{X}}
\newcommand{\BFV}{\bm{V}}
\newcommand{\BFS}{\bm{S}}

\newcommand{\BFehat}{\hat{\BFe}}
\newcommand{\BFGhat}{\widehat{\BFG}}
\newcommand{\BFShat}{\hat{\BFS}}
\newcommand{\BFshat}{\hat{\BFs}}

\newcommand{\sfH}{\mathsf{H}}
\newcommand{\sfI}{\mathsf{I}}

\newcommand{\sfQ}{\mathsf{Q}}
\newcommand{\sfU}{\mathsf{U}}
\newcommand{\sfV}{\mathsf{V}}

\newcommand{\sfHhat}{\widehat{\sfH}}

\newcommand{\mcA}{\mathcal{A}}

\newcommand{\mcE}{\mathcal{E}}
\newcommand{\mcF}{\mathcal{F}}
\newcommand{\mcG}{\mathcal{G}}

\newcommand{\mcK}{\mathcal{K}}

\newcommand{\mcN}{\mathcal{N}}
\newcommand{\mcO}{\mathcal{O}}

\newcommand{\mcS}{\mathcal{S}}

\newcommand{\mcU}{\mathcal{U}}

 \newcommand{\mcShat}{\hat{\mcS}}
  
\newcommand{\mbE}{\mathbb{E}}

\newcommand{\mbN}{\mathbb{N}}

\newcommand{\mbP}{\mathbb{P}}

\newcommand{\mbR}{\mathbb{R}}

\usepackage{color, colortbl}
\usepackage{bm}
\usepackage{mathrsfs}
\usepackage{booktabs}

\providecommand{\sfU}{\mathsf{U}}
\providecommand{\sfV}{\mathsf{V}}
\providecommand{\sfS}{\mathsf{S}}
\providecommand{\sfI}{\mathsf{I}}
\providecommand{\BFV}{\bm{V}}
\providecommand{\Ebar}{\bar{E}}
\providecommand{\Fbar}{\bar{F}}
\providecommand{\rhohat}{\hat{\rho}}
\providecommand{\sigmahat}{\hat{\sigma}}
\providecommand{\Dbar}{\bar{D}}
\providecommand{\Otilde}{\widetilde{\mathcal{O}}}
\providecommand{\ind}{\mathds{1}}
\providecommand{\Var}{\mathrm{Var}}
\providecommand{\mcU}{\mathscr{U}}
\providecommand{\mcA}{\mathscr{A}}
\providecommand{\mcK}{\mathscr{K}}
\providecommand{\mcG}{\mathscr{G}}
\providecommand{\mcN}{\mathscr{N}}

\providecommand{\mcE}{\mathscr{E}}

\providecommand{\kedF}{\kappa_{\mathrm{ed}F}}
\providecommand{\kQ}{\kappa_{\sfQ}}
\providecommand{\aQ}{\alpha_{\sfQ}}
\providecommand{\bQk}{\beta_{\sfQ,k}}
\providecommand{\bQ}{\beta_{\sfQ}}
\providecommand{\dQ}{
\bar{\beta}_{\sfQ}
}
\providecommand{\kH}{\kappa_{\sfH}}
\providecommand{\kLg}{\kappa_{\mathrm{L}g}}
\providecommand{\kLF}{\kappa_{\mathrm{L}F}}
\providecommand{\kg}{\kappa_{g}}
\providecommand{\kmf}{\kappa_{mf}}
\providecommand{\ksuc}{\kappa_{\mathrm{suc}}}
\providecommand{\ka}{\kappa_{a}}
\providecommand{\ta}{\theta_{a}}
\providecommand{\rF}{r_F}
\providecommand{\cF}{c_F}
\providecommand{\Ku}{K_{\mcU}}
\providecommand{\Ka}{K_{\mcA}}
\providecommand{\Kw}{K_{W}}
\providecommand{\SDelta}{S_{\Delta}}
\providecommand{\Teps}{T_{\epsilon}}
\providecommand{\Weps}{W_{\epsilon}}

\begin{document}

\title{Sample complexity and scalability of adaptive-sampling trust-region methods for zeroth-order stochastic optimization 
--- Part~II: a new pathwise analysis with randomized subspaces}
\author{S.~Shashaani\thanks{Email: sshasha2@ncsu.edu; North Carolina State University, NC, USA} \and C.~Cartis\thanks{University of Oxford, UK} \and K.~J.~Dzahini\thanks{Argonne National Laboratory, IL, USA}}

\date{}

\maketitle

\begin{abstract}
Stochastic trust-region methods are usually analyzed in expectation or with high probability: the guarantee is about the law of the algorithm, not about the run in hand. This paper analyzes the adaptive sampling stochastic trust-region method ASTRO-DF (which is the reuse variant ASTRA-R in Part I) in random subspaces along almost every sample path. The route has one probabilistic ingredient, the accuracy of the common-random-number estimates at the adaptive sample size, uniformly over the trust region; Part~I established it by chaining, here it is obtained at the stopping time by optional stopping and converted, by the first Borel--Cantelli lemma, into an accuracy that holds from a random index on, after which the deterministic trust-region argument applies. We prove almost sure convergence and almost sure iteration and sample complexity of orders $\tilde{\mathcal{O}}(d\epsilon^{-2})$ and $\tilde{\mathcal{O}}(d^2\epsilon^{-4})$, with a subspace dimension growing logarithmically; we show that at a fixed subspace dimension the eventual-accuracy argument yields   $\tilde{\mcO}(\epsilon^{-2-\mcO(1/p)})$ iterations and $\tilde{\mcO}(\epsilon^{-4-\mcO(1/p)})$ samples; we show that the random index after which the rates hold cannot be too large and has an explicit, rapidly decaying tail. Compared with the supermartingale route, the pathwise route needs no lower bound on the per-iteration probability of a good iteration, no coupling between the trust-region expansion factor and the sketch, and no moment condition on the summed squared radii.
\end{abstract}


\section{Introduction}\label{sec:intro}

We consider unconstrained simulation optimization problems of the form
\begin{equation}
    \min_{\BFx\in \mbR^d} f(\BFx) := \mbE[F(\BFx,\xi)] = \int_\Xi F(\BFx,\xi)\,dP(\xi),\label{eq:problem}
\end{equation}
where $f:\mbR^d\to\mbR$ is smooth, possibly nonconvex and bounded below, and $F:\mbR^d\times\Xi\to\mbR$ is observable only through a Monte Carlo oracle: a call at $\BFx$ returns $F(\BFx,\xi)=f(\BFx)+E(\BFx,\xi)$ for a fresh $\xi$, with no derivative. The estimator of $f(\BFx)$ from $n$ calls and its estimated variance are
\begin{equation}
\Fbar(\BFx,n)=\frac1n\sum_{j=1}^n F(\BFx,\xi_j),\qquad \sigmahat_F^2(\BFx,n)=\frac1n\sum_{j=1}^n\bigl(F(\BFx,\xi_j)-\Fbar(\BFx,n)\bigr)^2 .\label{eq:estimators}
\end{equation}
The oracle call is the expensive unit, and $W_\epsilon$ counts oracle calls. A \emph{replication} is one draw $\xi$ of the oracle's randomness, and a \emph{stream} is a sequence $\xi_1,\xi_2,\dots$ of independent replications. Evaluating several points on the same stream means evaluating $F(\cdot,\xi_j)$ at each of them with the same $\xi_j$, $j=1,\dots,n$; this is what common random numbers means here. A \emph{block} is a finite stretch of a stream, and a stream or block is \emph{fresh} if it is independent of everything drawn before it.

This is Part~II of a two-part study of one algorithm: the adaptive sampling stochastic trust-region method in random subspaces of Part~I~\cite{partI}, in the variant in which the trial point is evaluated on the replications already drawn to build the model. Part~I analyzed it through a supermartingale framework and obtained complexity bounds in expectation and almost surely along deterministic sequences, with a fixed subspace dimension. This part analyses the same algorithm along almost every sample path. The two parts share one probabilistic result, the uniform accuracy theorem of Part~I, and differ in everything built on top of it.

\paragraph{Two ways to evaluate efficiency in a stochastic algorithm}\label{sec:tworeadings}
A stochastic optimization algorithm run produces one sequence of iterates $\{\BFX_k\}$, one sequence of sample sizes $\{N_k\}$, and one first hitting time $\Teps:=\inf\{k:\|\nabla f(\BFX_k)\|\le\epsilon\}$, which are all random. There are two ways to make a statement about them.
 
The first, which is the common approach in the literature on stochastic trust-region methods~\cite{blanchet2019convergence,chen2018stochastic,jin2025sample,dzahini2024stochastic}, is to average over runs: $$\mbE[\Teps]\le c\,\epsilon^{-2},\qquad \text{or}\qquad \mbP(\Teps\le c\,\epsilon^{-2})\ge 1-e^{-c'\epsilon^{-2}}.$$ Such a statement is about the law of the algorithm. It is obtained by exhibiting a supermartingale in the iteration counter whose drift is negative as long as the algorithm has not stopped, and it requires that the good events of each iteration have conditional probability bounded below by a fixed constant, at every iteration. Its constants are deterministic. A run with $\Teps$ ten times the expected value is not excluded, only made unlikely.
 
The second is to fix the sample path and make a statement about it: for almost every $\omega$, there is a threshold $\epsilon_0(\omega)>0$ and a constant $R(\omega)<\infty$ such that $\Teps(\omega)\le R(\omega)\epsilon^{-2}$ for all $\epsilon\le\epsilon_0(\omega)$. Equivalently, along every deterministic sequence $\epsilon_\ell\downarrow 0$,
\begin{equation}
\mbP\bigl(T_{\epsilon_\ell}\,\epsilon_\ell^2>R\ \ \text{infinitely often}\bigr)=0 .\label{eq:asform}
\end{equation}
Such a statement is about the run. Its constant is random, and it says that on the run one did, beyond some random step (the ``burn-in'') the deterministic rate holds. In simulation optimization, where one typically gets a single optimization run with an expensive model, this is a guarantee worth having.
 
The two readings are obtained by different devices. Part~I obtains both: expected bounds, and almost-sure bounds of the form~\eqref{eq:asform} with an integrable constant, from a lower bound on the conditional probability of a good iteration holding at every iteration. This paper is about a different device for the second. 
The device is elementary. Every accuracy requirement of the deterministic trust-region theory (the model gradient is accurate to $\mcO(\Delta_k)$; the difference of estimation errors between incumbent and trial point is $\mcO(\Delta_k^2)$; the sketch is well aligned with the gradient) is an event $\mcE_k$ of iteration $k$. If its failure probability is summable, $\sum_k\mbP(\mcE_k^c)<\infty$, then by the first Borel--Cantelli lemma $\mcE_k$ holds for all $k\ge K(\omega)$ on almost every path, for a random, almost surely finite $K(\omega)$. From $K(\omega)$ onwards every hypothesis of the deterministic theory is in force, and the deterministic proof applies verbatim to the tail of the path. The first Borel--Cantelli lemma requires no independence, which is why common random numbers cause no difficulty on this route; and summability, not a lower bound at every iteration, is the only thing asked of the probabilities. What is lost is the expectation: $K(\omega)$ is finite and, as Section~\ref{sec:tails} shows, has finite moments, but the bound holds only beyond a random threshold $\epsilon_0(\omega)$, and no expected bound follows.
 
The cost of making failure probabilities summable is paid in two places. The sample size at iteration $k$ must be inflated by a slowly growing factor $\lambda_k$, so that the estimation failure probability is $e^{-c\lambda_k}$ and summable; $\lambda_k=\lambda_0(\log(k+2))^{1+\epsilon_\lambda}$ suffices, and it is the price of adaptive sampling known from sequential estimation~\cite{ChowRobbins1965,GhoshMukhopadhyaySen1997}. And the subspace dimension $p_k$ must grow so that the misalignment probability $e^{-cp_k}$ is summable; 
$p_k$ of order $\log k$, capped at $d$, suffices. Section~\ref{sec:fixedp} quantifies what happens if one refuses to pay the second price: at a fixed $p$ the argument of this paper gives $\epsilon^{-2-\mcO(1/p)}$, where Part~I’s per-iteration probability bound gives $\epsilon^{-2}$. 
 
It is fair to ask why one would take this route when the supermartingale route gives, for a fixed subspace dimension and correctly set constants, a stronger statement. The answer is in what each route asks of the algorithm rather than in what it delivers. 
The supermartingale route needs the good event of every iteration to have conditional probability above a fixed level; in particular the sketch must be well aligned with probability above one half at every iteration. The pathwise route needs only summability, so a sketch that is unreliable early and reliable late is admissible. 
The same feature admits components whose reliability improves with the iteration counter rather than being bounded below from the start, of which a growing subspace dimension is the instance treated here. And once the good events hold, the deterministic proof applies, so any result of the deterministic trust-region theory transfers to the stochastic setting by changing only the accuracy events; the supermartingale route re-derives a potential and a drift for each variant. Section~\ref{sec:tails} adds that the burn-in is not an unquantified random time: its tail is bounded by the tail of the series one sums in Borel--Cantelli, which is explicit, so finite-time, high-probability bounds are available from the same computation.

\subsection{Related work}\label{sec:related}
Stochastic trust-region methods with random models are analyzed through supermartingales in the STORM line of work~\cite{blanchet2019convergence,chen2018stochastic} and, in random subspaces, in~\cite{dzahini2024stochastic}; these analyses, and the high-probability ones of~\cite{cao2024first}, bound the number of iterations. The total sample complexity of STORM was obtained in~\cite{jin2025sample}, for sample sizes prescribed by the step size and a known variance bound. Adaptive sampling, in which the sample size is a stopping time driven by the observed variance, originates with ASTRO-DF~\cite{shashaani2018astro}; its pathwise complexity analysis~\cite{ha2025complexity} showed that common random numbers reduce the sample size per iteration from order $\Delta_k^{-4}$ to  $\Delta_k^{-2}$, and the sample complexity from order $\epsilon^{-6}$ to  $\epsilon^{-4}$. The same per-iteration reduction is obtained in~\cite{rinaldi2024stochastic,ding2025sequential} under a tail condition on the estimates, with deterministic sample sizes and no sample complexity bound. Part~I~\cite{partI} brought adaptive sampling and common random numbers into the supermartingale framework, including the reuse of replications at the trial point, with a fixed subspace dimension. Random-subspace model-based methods and their dimension dependence in the deterministic case are studied in~\cite{CRsubspace2021,cartis2026note}. The present paper is the pathwise counterpart of Part~I, with a subspace dimension that may grow, and it tracks the dimension dependence in the stochastic case.
 
\subsection{Contributions}\label{sec:contributions}
 
\emph{Pathwise analysis of the reuse algorithm.} Sections~\ref{sec:accuracy}--\ref{sec:complexity} give strong consistency and almost sure iteration and sample complexity, $\Teps=\Otilde(\epsilon^{-2})$ and $\Weps=\Otilde(\epsilon^{-4})$, for the algorithm of Part~I with fully reused common random numbers and a logarithmically growing subspace dimension. The analysis uses only three deterministic consequences of one probabilistic event (Section~\ref{sec:model}); it does not assume or prove full linearity of the model, and it does not use the accuracy of the model Hessian.
 
\emph{Concentration tailored to a stopping-time sample size.} The adaptive sample size $N_k$ is a stopping time for the within-iteration filtration. Section~\ref{sec:stopped} obtains the Bernstein tail of the averaged common-random-number increment at $N_k$ by optional stopping of the exponential supermartingale, with no statement at other sample sizes and no time-uniform maximal inequality. Part~I obtained the same tail by Ville's inequality, which is the stronger statement; we use the weaker one deliberately, because it uses exactly the property the algorithm supplies, and because it displays the mechanism behind the two steps discussed next. The chaining argument of Part~I, which makes the bound uniform over the trust region and is what permits reuse, is cited and not reproved; the one place where its proof changes is identified.
 
\emph{Relation to the pathwise analysis of ASTRO-DF.} This paper descends from~\cite{ha2025complexity}, which gave the pathwise analysis of ASTRO(-DF) in the full space and identified the sample-size scalings that common random numbers make possible. Two steps of its estimation-error argument are sharpened here, and both are instructive from the stopping-time viewpoint. Here we obtain the tail bound at the adaptive sample size by 
evaluating the exponential supermartingale at the stopping time itself (Section~\ref{sec:stopped}). And when incumbent and candidate share replications, the candidate is a function of the replications being averaged; a pointwise inequality applies to it when it is evaluated on a random number stream drawn after the model is built, with the pairing kept within that stream, 
while full reuse calls for the uniform bound of Section~\ref{sec:uniform}. Section~\ref{sec:repairs} sets out both points and the three sampling schemes they distinguish.
 
\emph{Fixed versus growing subspace dimension.} Section~\ref{sec:fixedp} proves that with a fixed $p$ and misalignment probability $\dQ$ below an explicit threshold, the radius satisfies $\Delta_k\ge b_\epsilon k^{-\varrho}$ eventually with $\varrho=\mcO(1/\log(1/\dQ))$, which gives almost sure liminf convergence when $\varrho\le\frac12$ and $\Teps=\mcO(\epsilon^{-2/(1-2\varrho)})$. The exponent tends to $2$ as $p$ grows, 
and $p_k=\mcO(\log k)$ is the threshold at which the exponent is exactly~$2$. The growth is what the eventual-accuracy argument needs to reach exponent $2$ without Part~I’s per-iteration probability bound. 
 
\emph{Dimension dependence.} Section~\ref{sec:dimension} discusses all statements about $d$, $p$ and the sketch scale. Rescaling the sketch is shown to be a reparametrization of the algorithm (Proposition~\ref{prop:reparam}), so that the improvement of Cartis and Roberts~\cite{cartis2026note} is the improvement of one parametrization over another; under their convention orthonormal columns give $\Otilde(d\epsilon^{-2})$ iterations, the sketch of~\cite{dzahini2024stochastic} scaled by $\sqrt{d/p}$ gives $\Otilde(d^2p^{-1}\epsilon^{-2})$, and a scale balancing model bias against estimation error gives $\Otilde(dp^{-1/2}\epsilon^{-2})$. 
The price of reuse in our bounds is a factor $p$ in the sample size per iteration: the accuracy must hold at every point of the $p$-dimensional trust region at once, and the nets that cover it have a number of points exponential in $p$. With $p_k=\mcO(\log k)$ this factor is logarithmic; in the full space, $p=d$, it is a factor $d$. This scale, $\kQ\asymp p^{-1/4}$, is the best over the sketch scale in our bounds and gains a factor $\sqrt p$ over orthonormal columns in both the iteration and the sample bound. It requires the tolerance on the estimation error to shrink with the subspace dimension, $\kedF\asymp p^{-1/2}$, which costs a factor $\sqrt p$ in the sample size at a given radius.

During the writing of this work, a recent survey~\cite{berahas2026stochastic} appeared
that formalizes common random numbers as a \emph{multi-point} stochastic zeroth-order oracle. The ASTRO-DF \cite{ha2025complexity} and ASTRA \cite{partI} sample complexity bounds of order $\mathcal{O}(\epsilon^{-4})$ are achieved in \cite{berahas2026stochastic} in expectation, for a fixed-step randomized finite-difference method. Here and in \cite{partI}, we address the more analytically-challenging and practically-useful framework of stochastic trust region methods with parameter-free adaptive sampling, with fresh, or even re-used, objective estimates. The dimension dependence of our bounds here matches the ones in \cite{berahas2026stochastic}.

\subsection{Relation to Part~I}\label{sec:partI}
 
Both parts analyze Algorithm~\ref{alg:main} below. Table~\ref{tab:parts} records what each requires and delivers. Part~I asks that the good event of each iteration have conditional probability above a fixed level at every iteration, and needs the misalignment probability of the sketch below $\frac12$; from this it obtains, through the supermartingale on the radius exponent, constants with finite expectation and bounds in expectation. Part~II asks that the failure probabilities be summable and obtains, through Borel--Cantelli, that the good events hold eventually; the constants are random (integrable when $\delta>1$ in~\eqref{eq:pk}, Remark~\ref{rem:asScope}), the subspace dimension must grow, and nothing is obtained in expectation. 

\vspace{1em}
\begin{table}[h]\centering\small
\caption{Requirements and outputs of the two routes. Both use Algorithm~\ref{alg:main} and the uniform accuracy bound of Theorem~\ref{thm:uniform} (Part~I, Theorem~4.2). In Part~I, $\bQ$ is the misalignment probability; here $\bQk$ is the alignment probability (Assumption~\ref{assum:well-aligned}). 
}\label{tab:parts}
\begin{tabular}{l p{6.2cm} p{6.5cm}}
\toprule
& Part~I (supermartingale) & Part~II (pathwise)\\
\midrule
within an iteration & Ville's inequality & optional stopping (Lemma~\ref{lem:stopped})\\
across iterations & count of favorable events, each of conditional probability $\ge\pi_\star>\frac12$ & first Borel--Cantelli; summable failure probabilities\\
subspace dimension & fixed $p$, with $\bQ<1-\hat\gamma_r/\pi_{\min}$ & $p_k=\mcO(\log k)$ nondecreasing; fixed $p$ at a cost (\S\ref{sec:fixedp})\\
parameter coupling & $\pi_\star\gamma^{-2r}+(1-\pi_\star)\gamma^{2r}<1$ & none\\
inflation $\lambda_k\to\infty$ & optional; a constant inflation needs thresholds on $\lambda_0$ set by unknown constants & required\\
constants $\ka,\kg$ & fixed & may vary with $k$\\
iteration bound & $\mcO(\epsilon^{-2})$ a.s.\ and in expectation & $\Otilde(\epsilon^{-2})$ a.s.\\
constants & deterministic, but through the index beyond which its probability bounds hold, which unknown problem constants fix & random; integrable when $p_k$ grows a little faster than the minimum (\eqref{eq:pk} with $\delta>1$, Remark~\ref{rem:asScope})\\
range in $\epsilon$ & expectation only below a threshold exponentially small in that index & every $\epsilon$ below a random threshold $\epsilon_0(\omega)$\\
sample bound & $\Otilde(\epsilon^{-4})$ a.s.; in expectation under a moment condition on $\SDelta$ & $\Otilde(\epsilon^{-4})$ a.s.\\
tail of the burn-in & Markov from the expectation, order $t^{-1}$ & explicit, superpolynomial (\S\ref{sec:tails})\\
noise moments & fourth & fourth\\
dimension dependence & not treated & \S\ref{sec:dimension}\\
\bottomrule
\end{tabular}
\end{table}

Both routes use Assumptions~\ref{assum:smooth}--\ref{assum:cauchy}, \ref{assum:well-aligned}(b) and, for the sample complexity bound, \ref{assum:fourth}. The supermartingale route of Part~I asks in addition that the conditional probability of a good iteration be bounded below, at every iteration beyond a deterministic index, by a level exceeding one half (its Assumption~2.1(i)), which with its sketch is a misalignment probability below $1-\hat\gamma_r/\pi_{\min}<\frac12$ that a fixed $p$ must meet; that $\pi_\star\gamma^{-2r}+(1-\pi_\star)\gamma^{2r}<1$ (its Assumption~2.2), which couples the expansion factor of the trust region to the reliability of the sketch; and, for its expected sample bound under the logarithmic inflation, a moment condition on $\SDelta$ beyond $\mbE[\SDelta]<\infty$ that it assumes rather than derives (its Theorem~6.4). The pathwise route asks none of these when $p_k$ grows: alignment need only fail summably often, nothing couples $\gamma$ to the sketch, and no moment of $\SDelta$ enters; the only probabilistic input is the unconditional summability $\sum_k\mbP(\mcU_k^c)<\infty$. It asks instead for $\lambda_k\to\infty$ and for $p_k\to\infty$, the latter at the price of Section~\ref{sec:fixedp} if refused, and it gives up 
statements in expectation. It also admits constants that vary with $k$, which the framework of Part~I does not, and this is how a growing subspace dimension enters at all.

\subsection{Notation}\label{sec:notation}

Deterministic scalars are lower case ($\kappa$, $c$, $\theta$, $\lambda$, $r_F$), random scalars are upper case ($\Delta_k$, $N_k$, $\Teps$, the burn-in indices $\Ku,\Ka,\Kw$ and the constants $R,\SDelta$ below are random). Vectors are bold, random vectors bold upper case ($\BFX_k$), matrices sans serif ($\sfQ_k$), sets and $\sigma$-fields script. Hatted symbols live in the subspace. Constants depending on the subspace dimension are written as functions of $p$, as $\kg(p)$, and evaluated at $p_k$ when used at iteration $k$. We write $a\wedge b=\min\{a,b\}$, $a\vee b=\max\{a,b\}$, $\{\mcE_k\ \text{i.o.}\}=\limsup_k\mcE_k$, and $\Otilde(\cdot)$ for an order up to a factor polylogarithmic in the argument. All random quantities live on one probability space $(\Omega,\mathscr F,\mbP)$.

\section{The algorithm}\label{sec:algorithm}

\subsection{Interpolation on a coordinate design in a subspace}\label{sec:design}

At iteration $k$ the method holds an incumbent $\BFX_k^0\in\real^d$, a radius $\Delta_k>0$, and draws a sketch $\sfQ_k\in\real^{d\times p_k}$. It searches the affine subspace $\{\BFX_k^0+\sfQ_k\BFshat:\ \BFshat\in\real^{p_k}\}$ within $\|\BFshat\|\le\Delta_k$, and writes $\fhat_k(\BFshat):=f(\BFX_k^0+\sfQ_k\BFshat)$, so that $\nabla\fhat_k(\BFshat)=\sfQ_k^\intercal\nabla f(\BFX_k^0+\sfQ_k\BFshat)$. The design is the coordinate design \emph{of the subspace},
\begin{equation}
    \mcShat_k=\{\bm{0}_{p_k},\ \pm\Delta_k\BFehat^1,\dots,\pm\Delta_k\BFehat^{p_k}\},\qquad \BFX_k^i:=\BFX_k^0+\sfQ_k\BFShat_k^i,\quad i=0,1,\dots,2p_k,\label{eq:interpolation-set}
\end{equation}
with $\BFShat_k^i=\Delta_k\BFehat^i$ and $\BFShat_k^{p_k+i}=-\Delta_k\BFehat^i$ for $i\in\{1,\dots,p_k\}$. In the ambient space the design points lie within $\|\sfQ_k\|\Delta_k$ of $\BFX_k^0$. The full-space method is $\sfQ_k=\sfI_d$, $p_k=d$.

The model is quadratic with diagonal Hessian, and it is determined by the $2p_k+1$ estimates at the design points: with $\Fbar_k^i:=\Fbar(\BFX_k^i,N_k)$,
\begin{equation}
    [\BFGhat_k]_i=\frac{\Fbar_k^i-\Fbar_k^{p_k+i}}{2\Delta_k},\qquad
    [\sfHhat_k]_{ii}=\kH\wedge\Bigl(\bigl(-\kH\bigr)\vee\frac{\Fbar_k^i-2\Fbar_k^0+\Fbar_k^{p_k+i}}{\Delta_k^2}\Bigr),\qquad i=1,\dots,p_k,\label{eq:GH}
\end{equation}
and
\begin{equation}
    \Mhat_k(\BFshat)=\Fbar_k^0+\BFshat^\intercal\BFGhat_k+\tfrac12\BFshat^\intercal\sfHhat_k\BFshat .\label{eq:model}
\end{equation}
The gradient is a central finite difference along the columns of $\sfQ_k$, 
and the diagonal Hessian is a second difference truncated to $[-\kH,\kH]$, that is, each diagonal entry is replaced by the nearest point of that interval, for a constant $\kH\ge1$, so that $\|\sfHhat_k\|\le\kH$ surely. The truncation is a safeguard for inaccurate iterations, on which a second difference of noisy estimates, an estimation error divided by $\Delta_k^2$, can be arbitrarily large. On the good event of Section~\ref{sec:accuracy} it changes nothing once $\kH\ge\kLg\kQ^2+2\kedF$ (Lemma~\ref{lem:hessian}), so the model is eventually the plain interpolation model; 
it replaces the bounded-Hessian assumption of~\cite{ha2025complexity}. The second difference is taken in subspace coordinates, divided by $\Delta_k^2$ and not by $\|\sfQ_k\BFehat^i\|^2\Delta_k^2$, because $\sfHhat_k$ approximates $\sfQ_k^\intercal\nabla^2f\,\sfQ_k$ and not $\nabla^2 f$; the other normalization makes the estimation part of the Hessian bound grow as the sketch is scaled down and spoils the accounting of Section~\ref{sec:dimension}.

\subsection{One stream, one rule, full reuse}\label{sec:onestream}

All estimates at iteration $k$ use the same replications $\xi_{k,1},\xi_{k,2},\dots$, and sampling stops at
\begin{equation}
    N_k=\min\Bigl\{n\in\mbN:\ \frac{\sigma_0\vee\sigmahat^{\max}_F(k,n)}{\sqrt n}\le\kappa_s\frac{\Delta_k}{\sqrt{\lambda_k}}\Bigr\},\qquad \sigmahat^{\max}_F(k,n):=\max_{0\le i\le 2p_k}\sigmahat_F(\BFX_k^i,n),\label{eq:stopping}
\end{equation}
with $\lambda_k=\lambda_0(\log(k+2))^{1+\epsilon_\lambda}$, $\epsilon_\lambda\in(0,1)$, and constants $\sigma_0,\kappa_s>0$. The rule compares the largest observed standard error over the design with a threshold set by the algorithm's own radius; it requires no knowledge of the noise level, no Lipschitz constant and no target accuracy. Since $\sigma_0\vee\sigmahat^{\max}_F\ge\sigma_0$, any $n$ meeting the criterion satisfies
\begin{equation}
    N_k\ \ge\ n_k:=\Bigl\lceil\frac{\sigma_0^2\lambda_k}{\kappa_s^2\Delta_k^2}\Bigr\rceil ,\label{eq:floor}
\end{equation}
a floor that is determined before the first replication of the iteration is drawn. The floor is the only property of the rule that the accuracy analysis uses; what the rule adds to the floor is efficiency, see Remark~\ref{rem:ruleBuys}.

The trial point is evaluated on the same replications: $\Fbar_k^{\text{s}}:=\Fbar(\BFX_k^{\text{s}},N_k)$ with $\BFX_k^{\text{s}}=\BFX_k^0+\sfQ_k\BFShat_k^*$, and no further sampling takes place in the iteration. This is what an implementation does, and it is the cheapest of the three schemes of Remark~\ref{rem:no-crn} at equal sample sizes: $(2p_k+2)N_k$ calls per iteration. It is also what the existing supermartingale analyses (prior to Part~I) exclude, since they require the trial estimates to be accurate conditionally on the $\sigma$-field generated by the model, which the trial point depends on~\cite{chen2018stochastic,dzahini2024stochastic}, and what the pathwise analysis of~\cite{ha2025complexity} treated with a pointwise inequality, which Section~\ref{sec:repairs} revisits.

\begin{algorithm}[htp]\small
\caption{ASTRO-DF in random subspaces with fully reused CRN}
\label{alg:main}
\begin{algorithmic}[1]
\REQUIRE $\BFx_0\in\real^d$; radii $\Delta_{\max}>\Delta_0>0$; thresholds $\eta\in(0,1)$, $\mu>0$; expansion and contraction constants $\gamma_{inc}>1>\gamma_{dec}>0$; Hessian clip $\kH\ge1$; inflation $\lambda_k=\lambda_0(\log(k+2))^{1+\epsilon_\lambda}$; sampling constants $\sigma_0,\kappa_s>0$; subspace dimensions $\{p_k\}$.
\FOR{$k=0,1,2,\dots$}
\STATE \label{ASTRO:sketch}\textit{Sketch:} draw $\sfQ_k\in\real^{d\times p_k}$ independently of the past and of all replications.
\STATE \label{ASTRO:Model}\textit{Model:} evaluate the design points~\eqref{eq:interpolation-set} on one stream $\xi_{k,1},\xi_{k,2},\dots$; stop at $N_k$ from~\eqref{eq:stopping}; form $\BFGhat_k$, $\sfHhat_k$ by~\eqref{eq:GH} and $\Mhat_k$ by~\eqref{eq:model}.
\STATE \label{ASTRO:TRsubprob}\textit{Subproblem:} compute $\BFShat_k^*\approx\argmin_{\|\BFshat\|\le\Delta_k}\Mhat_k(\BFshat)$ satisfying Assumption~\ref{assum:cauchy}; set $\BFX_k^{\text{s}}=\BFX_k^0+\sfQ_k\BFShat_k^*$.
\STATE \label{ASTRO:evaluate}\textit{Trial evaluation:} evaluate $\BFX_k^{\text{s}}$ on $\xi_{k,1},\dots,\xi_{k,N_k}$, set $\Fbar_k^{\text{s}}:=\Fbar(\BFX_k^{\text{s}},N_k)$ and
\[\rhohat_k=\frac{\Fbar_k^0-\Fbar_k^{\text{s}}}{\Mhat_k(\bm{0})-\Mhat_k(\BFShat_k^*)}.\]
\STATE \label{ASTRO:ratio}\textit{Update:} $(\BFX_{k+1}^0,\Delta_{k+1})=(\BFX_k^{\text{s}},\ \gamma_{inc}\Delta_k\wedge\Delta_{\max})$ if $\rhohat_k\ge\eta$ and $\mu\|\BFGhat_k\|\ge\Delta_k$ (success); otherwise $(\BFX_k^0,\ \gamma_{dec}\Delta_k)$ (failure).
\ENDFOR
\end{algorithmic}
\end{algorithm}

\subsection{Filtrations}\label{sec:filtrations}

Let $\mcF_{k-1}$ be the $\sigma$-field generated by everything determined by the end of iteration $k-1$, so that $\BFX_k^0$ and $\Delta_k$ are $\mcF_{k-1}$-measurable, and set
\begin{equation}
    \mcF^{\sfQ}_{k-1}:=\mcF_{k-1}\vee\sigma(\sfQ_k),\qquad \mcG_{k,n}:=\mcF^{\sfQ}_{k-1}\vee\sigma(\xi_{k,1},\dots,\xi_{k,n}),\qquad \mcG_{k,0}=\mcF^{\sfQ}_{k-1},\label{eq:filtrations}
\end{equation}
with $\mcF_k=\mcG_{k,N_k}$. The sketch, the design points and the alignment event of Assumption~\ref{assum:well-aligned} are $\mcF^{\sfQ}_{k-1}$-measurable; the model, the step, the trial point and $\rhohat_k$ are $\mcG_{k,N_k}$-measurable. No $\sigma$-field generated by the model appears anywhere in the analysis, and this is not an accident: 
under reuse the summands of the trial estimate are not centered conditionally on the replications already drawn, so no accuracy statement conditional on the model is available, 
and the analysis must be organized so as never to need one.

\subsection{The primitive object}\label{sec:primitive}

Common random numbers act on the paired increment
\begin{equation}
    D(\BFx,\BFy;\xi):=\bigl(F(\BFy,\xi)-f(\BFy)\bigr)-\bigl(F(\BFx,\xi)-f(\BFx)\bigr),\qquad \mbE_\xi[D]=0,\label{eq:D}
\end{equation}
and for $\|\BFshat\|\le\Delta_k$ we write
\begin{equation}
    \Dbar_k(\BFshat,n):=\frac1n\sum_{j=1}^n D\bigl(\BFX_k^0,\ \BFX_k^0+\sfQ_k\BFshat;\ \xi_{k,j}\bigr),\label{eq:Dbar}
\end{equation}
so that $\Dbar_k(\BFShat_k^i,N_k)=\Ebar_k^i-\Ebar_k^0$ at a design point and $\Dbar_k(\BFShat_k^*,N_k)=\Ebar_k^{\text{s}}-\Ebar_k^0$ at the trial point, where $\Ebar_k^i:=\Fbar_k^i-f(\BFX_k^i)$. The identity
\begin{equation}
    D(\BFx,\BFy;\xi)-D(\BFx,\BFy';\xi)=D(\BFy',\BFy;\xi)\label{eq:telescope}
\end{equation}
says that increments of $\BFshat\mapsto\Dbar_k(\BFshat,n)$ are themselves paired increments; it is what makes the process chainable.

\section{Assumptions}\label{sec:assumptions}

We use the final assumptions of Part~I, and we state them in full so that this part can be read alone.

\begin{assumption}[smooth $f$]\label{assum:smooth}
$f$ is continuously differentiable on $\real^d$ and bounded below by $f^*:=\inf_{\BFx}f(\BFx)>-\infty$, and $\nabla f$ is $\kLg$-Lipschitz.
\end{assumption}

\begin{assumption}[finite variance]\label{assum:variance}
For all $\BFx\in\real^d$, $\sigma_F^2(\BFx):=\Var(F(\BFx,\xi))\le v_F<\infty$.
\end{assumption}

\begin{assumption}[Bernstein condition on the paired increment]\label{assum:bernstein}
There are constants $b_D,\sigma_D>0$ such that for all $\BFx,\BFy\in\real^d$ with $\|\BFy-\BFx\|\le\kQ\Delta_{\max}$ and every integer $m\ge2$,
\begin{equation}
    \mbE\bigl[|D(\BFx,\BFy;\xi)|^m\bigr]\le\frac{m!}{2}\bigl(b_D\|\BFy-\BFx\|\bigr)^{m-2}\bigl(\sigma_D\|\BFy-\BFx\|\bigr)^2 .\label{eq:bernsteinD}
\end{equation}
\end{assumption}

   \begin{assumption}[continuous sample paths]\label{assum:lip-F}
   The sample paths $\BFx\mapsto F(\BFx,\xi)$ are continuous almost surely.
   \end{assumption}

\begin{remark}[a single sufficient condition, and what changed since~\cite{ha2025complexity}]\label{rem:singleSufficient}
Assumptions~\ref{assum:bernstein} and~\ref{assum:lip-F} are both implied by Lipschitz sample paths, that is, $|F(\BFx_1,\xi)-F(\BFx_2,\xi)|\le\kLF(\xi)\|\BFx_1-\BFx_2\|$, (Lipschitz paths are continuous) with a light-tailed modulus, $\mbE[e^{\theta_0\kLF(\xi)}]<\infty$ for some $\theta_0>0$: writing $$W:=\kLF(\xi)+\mbE[\kLF(\xi)],$$ Jensen's inequality gives $|D(\BFx,\BFy;\xi)|\le W\|\BFy-\BFx\|$, and $e^x\ge x^m/m!$ gives $$\mbE[W^m]\le\frac{m!}{2}(\theta_0^{-1})^{m-2}\cdot2\theta_0^{-2}\mbE[e^{\theta_0W}],$$ which is~\eqref{eq:bernsteinD} with $b_D=\theta_0^{-1}$ and $\sigma_D^2=2\theta_0^{-2}\mbE[e^{\theta_0 W}]$. 
Assumption~4 of~\cite{ha2025complexity} carried the same information through the constant $\kappa_{\mathrm{L}m}:=\max_{m\ge2}(\mbE|\kLF|^m)^{1/m}$; since $\sup_m\|X\|_m=\|X\|_\infty$ for any random variable, that formulation amounts to a bounded modulus, and the Bernstein form~\eqref{eq:bernsteinD}, adopted in that paper as well, is the one that admits an unbounded light-tailed modulus, as arises in simulation models with heavy-tailed service or interarrival times. 
\end{remark}

\begin{assumption}[Cauchy decrease in the subproblem]\label{assum:cauchy}
For some $\kappa_{\mathrm{fcd}}\in(0,1]$ and all $k$,
$$\Mhat_k(\bm{0})-\Mhat_k(\BFShat_k^*)\ge\tfrac12\kappa_{\mathrm{fcd}}\|\BFGhat_k\|\bigl(\|\BFGhat_k\|/(\|\sfHhat_k\|\vee1)\wedge\Delta_k\bigr).$$
\end{assumption}

\begin{assumption}[well-aligned sketches]\label{assum:well-aligned}
The matrices $\{\sfQ_k\}$ are drawn independently of $\{\xi_{k,j}\}$ and satisfy, for constants $\aQ\in(0,1]$, $\kQ<\infty$ and a deterministic sequence $\bQk\in[0,1]$:
\begin{enumerate}
\item[(a)] $\mbP\bigl(\|\sfQ_k^\intercal\BFV\|\ge\aQ\|\BFV\|\ \big|\ \mcF_{k-1}\bigr)\ge\bQk$ almost surely, for every $\mcF_{k-1}$-measurable random vector $\BFV\in\real^d$;
\item[(b)] $\|\sfQ_k\|\le\kQ$ almost surely.
\end{enumerate}
We write $$\mcA_k:=\{\|\sfQ_k^\intercal\nabla f(\BFX_k^0)\|\ge\aQ\|\nabla f(\BFX_k^0)\|\}\in\mcF^{\sfQ}_{k-1}$$ and $\dQ:=\sup_k(1-\bQk)$ for the misalignment probability.
\end{assumption}

Part~(a) is stated conditionally, and for $\mcF_{k-1}$-measurable $\BFV$, because it is applied to $\nabla f(\BFX_k^0)$, which depends on $\sfQ_0,\dots,\sfQ_{k-1}$; the conditional form follows from the unconditional one for fixed $\BFv$ because $\sfQ_k$ is independent of $\mcF_{k-1}$. Note that $\bQk$ is a \emph{success} probability, whereas Part~I and~\cite{dzahini2024stochastic} parametrize by the failure probability. The pathwise route uses (a) in one of two regimes: $\sum_k(1-\bQk)<\infty$, achieved by a growing subspace dimension (Sections~\ref{sec:conv}--\ref{sec:complexity}), or $\dQ<1$ fixed (Section~\ref{sec:fixedp}). Which ensembles satisfy (a) and (b), and with which constants, is the subject of Section~\ref{sec:sketch}.

\begin{assumption}[fourth moment; used for the sample complexity bound only]\label{assum:fourth}
There is $c_4<\infty$ with $$\sup_{\BFx\in\real^d}\mbE\bigl[|F(\BFx,\xi)-f(\BFx)|^4\bigr]\le c_4.$$
\end{assumption}

\begin{remark}\label{rem:momentScope}
Assumption~\ref{assum:fourth} is used only in Lemma~\ref{lem:Nk}, to bound the realized sample size from above, and hence only in the sample complexity Theorem~\ref{thm:work}. Convergence and iteration complexity do not need it. It is confined because it can fail exactly where common random numbers are most useful: if $F(\BFx,\xi)=f(\BFx)+Z+R(\BFx,\xi)$ with $Z$ a common additive noise of finite variance but no fourth moment and $R$ having light-tailed increments proportional to $\|\BFy-\BFx\|$, then $Z$ cancels in~\eqref{eq:D}, Assumptions~\ref{assum:variance}--\ref{assum:lip-F} hold, and $\mbE|F-f|^4=\infty$. For such noise the sample size $N_k$, which is driven by the empirical variance of $F$ and therefore by $Z$, cannot be bounded by the paired increments, and a bound on the sample complexity does require control of $\sigmahat_F^2$ that only a moment condition on $F$ itself supplies. The fourth moment is sufficient for the deviation bound of Lemma~\ref{lem:Nk}; we do not claim it is necessary. The summability used in Lemma~\ref{lem:Nk} comes from $\sum_k\Delta_k^2<\infty$, and a $(2+\delta)$-th moment would require $\sum_k\Delta_k^\delta<\infty$, which is not available.
\end{remark}

\subsection{Sketch ensembles}\label{sec:sketch}

The analysis uses the sketch only through Assumption~\ref{assum:well-aligned}. The design is the coordinate design of the subspace, that is, the points $\BFX_k^0\pm\Delta_k\sfQ_k\BFehat^i$ along the columns of $\sfQ_k$, and nothing is assumed about those columns beyond (a) and (b). Two consequences deserve to be made explicit. First, (a), (b) and the constants $\kg,\ka$ below are invariant under rotations within the subspace (the iterates are not, since the coordinate design turns with the basis): replacing $\sfQ_k$ by $\sfQ_k\sfV_k$ for an $\mcF^{\sfQ}_{k-1}$-measurable orthogonal $\sfV_k\in\real^{p_k\times p_k}$ composes $\fhat_k$ with a rotation and changes neither (a) nor (b); the factor $\sqrt{p}$ in the gradient error (Lemma~\ref{lem:center}) comes from the $2p$ rows of the design matrix in subspace coordinates and is rotation invariant. Second, sketches with orthonormal columns are a special case, and the convenient one.

\begin{proposition}[orthonormal columns]\label{prop:Haar}
Let $\sfU_k\in\real^{d\times p_k}$ consist of the first $p_k$ columns of a matrix drawn from the Haar measure on the orthogonal group $\mathbb O(d)$, independently of $\mcF_{k-1}$ and of all replications, and let $\sfQ_k:=c\,\sfU_k$ for a constant $c>0$. Then $\|\sfQ_k\|=c$ surely, so (b) holds with $\kQ=c$. Moreover, for every $\alpha\in(0,1)$ and every $\mcF_{k-1}$-measurable $\BFV$, almost surely,
\begin{equation}
\mbP\Bigl(\|\sfQ_k^\intercal\BFV\|\ge c\sqrt{(1-\alpha)p_k/d}\ \|\BFV\|\ \Big|\ \mcF_{k-1}\Bigr)\ \ge\
\begin{cases}1-e^{-p_k\alpha^2/4},&p_k<d,\\ 1,&p_k=d.\end{cases}\label{eq:HaarAlign}
\end{equation}
Consequently, if $p_k\ge p_0$ for all $k$ and $c\le\sqrt{d/((1-\alpha)p_0)}$, then (a) holds with $\aQ=c\sqrt{(1-\alpha)p_0/d}\in(0,1]$, with $1-\bQk=e^{-p_k\alpha^2/4}$ when $p_k<d$ and $1-\bQk=0$ when $p_k=d$.
\end{proposition}
\begin{proof}
Orthonormality of the columns gives $\|\sfQ_k\BFshat\|=c\|\BFshat\|$ for every $\BFshat$, hence $\|\sfQ_k\|=c$. Fix $\alpha$ and $\BFV$. On $\{\BFV=\bm0\}$ the event in~\eqref{eq:HaarAlign} holds trivially. If $p_k=d$, $\sfU_k$ is orthogonal, so $\|\sfQ_k^\intercal\BFV\|=c\|\BFV\|$ and the event holds surely. Let $p_k<d$. Since $\sfU_k$ is independent of $\mcF_{k-1}$ and $\BFV$ is $\mcF_{k-1}$-measurable, the conditional probability in~\eqref{eq:HaarAlign} equals $h(\BFV)$ with $$h(\BFv):=\mbP\bigl(\|\sfU_k^\intercal\BFv\|^2\ge(1-\alpha)(p_k/d)\|\BFv\|^2\bigr)$$ for deterministic $\BFv$. For $\BFv\ne\bm0$, rotation invariance of the Haar measure shows that $\|\sfU_k^\intercal\BFv\|^2/\|\BFv\|^2$ has the law of the squared norm of the first $p_k$ coordinates of a uniformly distributed point on the unit sphere of $\real^d$, and~\cite[Lemma~2.2]{DasguptaGupta2003JLT} gives
\[
1-h(\BFv)\le\exp\Bigl(\frac{p_k}{2}\bigl(\alpha+\log(1-\alpha)\bigr)\Bigr)\le e^{-p_k\alpha^2/4},
\]
because $\log(1-\alpha)\le-\alpha-\alpha^2/2$. This proves~\eqref{eq:HaarAlign}. Finally, $p_k\ge p_0$, so the event in~\eqref{eq:HaarAlign} is contained in $\{\|\sfQ_k^\intercal\BFV\|\ge c\sqrt{(1-\alpha)p_0/d}\,\|\BFV\|\}$, which is (a) with the stated constants.
\end{proof}

We now make three remarks:
\begin{itemize}
    \item[(i)] Proposition~\ref{prop:Haar} makes (b) hold \emph{surely}; for a Gaussian Johnson--Lindenstrauss transform $\|\sfQ_k\|$ is unbounded and (b) holds only with high probability, so the orthonormal construction is the one for which the analysis applies without a further union bound. Hashing transforms satisfy (b) surely as well~\cite{cartis2026note}. 
    \item[(ii)] The design points then lie at exact distance $c\Delta_k$ from $\BFX_k^0$, and $\fhat_k$ is $f$ restricted to a $p_k$-plane and dilated by $c$. 
    \item[(iii)] The ensemble-intrinsic quantity is $\kQ^2/\aQ^2\asymp d/p$: for orthonormal columns $\kQ=c$ and $\aQ\asymp c\sqrt{p/d}$, and for a Gaussian or hashing transform $\sfS^\intercal\in\real^{d\times p}$ 
    one has $\kQ\asymp c\sqrt{d/p}$ and $\aQ\asymp c$~\cite{cartis2026note}. 
    The two ensembles enter the analysis through this quantity alone, and the remaining freedom is the scale $c$. The next proposition says what that freedom is.
\end{itemize}

\begin{proposition}[rescaling the sketch is a reparametrization]\label{prop:reparam}
Fix $c>0$. Run Algorithm~\ref{alg:main} once with sketches $\{\sfQ_k\}$ and parameters $(\Delta_0,\Delta_{\max},\mu,\kappa_s,\kH)$, and once with sketches $\{c\,\sfQ_k\}$, the same replications, and parameters $$(\Delta_0/c,\ \Delta_{\max}/c,\ \mu/c^2,\ c\,\kappa_s,\ c^2\kH),$$ and a subproblem solver that commutes with dilation of the subspace coordinates (the exact minimizer, the Cauchy point, dogleg and truncated CG do). The two runs produce the same iterates $\BFX_k^0$, trial points $\BFX_k^{\text{s}}$, sample sizes $N_k$, ratios $\rhohat_k$ and success indicators for every $k$, and their radii satisfy $\Delta_k'=\Delta_k/c$. Moreover the event $\mcU_k$ of~\eqref{eq:Uk} for the second run with accuracy constant $c^2\kedF$ is the event $\mcU_k$ of the first run with constant $\kedF$.
\end{proposition}

\begin{proof}
Induction on $k$, with primes for the second run and $\Delta_0'=\Delta_0/c$. With $\Delta_k'=\Delta_k/c$ the design points coincide, $\BFX_k^0\pm\Delta_k'(c\sfQ_k)\BFehat^i=\BFX_k^0\pm\Delta_k\sfQ_k\BFehat^i$, so the estimates coincide; the rule~\eqref{eq:stopping} reads $$\sigmahat/\sqrt n\le c\kappa_s\Delta_k'/\sqrt{\lambda_k}=\kappa_s\Delta_k/\sqrt{\lambda_k},$$ so $N_k'=N_k$. By~\eqref{eq:GH}, $\BFGhat_k'=c\,\BFGhat_k$ and the unclipped second differences scale by $c^2$, which the clip level $c^2\kH$ preserves; hence $\Mhat_k'(\BFshat')=\Mhat_k(c\BFshat')$, minimizing over $\|\BFshat'\|\le\Delta_k'$ is minimizing $\Mhat_k$ over $\|\BFshat\|\le\Delta_k$ with $\BFshat=c\BFshat'$, the trial points coincide, and Assumption~\ref{assum:cauchy} holds in the second run with $\kappa_H$ read as $c^2\kH\vee1$. The ratios coincide, and $(\mu/c^2)\|\BFGhat_k'\|\ge\Delta_k'$ is $\mu\|\BFGhat_k\|\ge\Delta_k$. Both runs take the same decision and $\Delta_{k+1}'=\Delta_{k+1}/c$ since $\Delta_{\max}'=\Delta_{\max}/c$. Finally $$\sup_{\|\BFshat'\|\le\Delta_k'}|\Dbar_k'(\BFshat',N_k)|=\sup_{\|\BFshat\|\le\Delta_k}|\Dbar_k(\BFshat,N_k)|$$ while $c^2\kedF\Delta_k'^2=\kedF\Delta_k^2$.
\end{proof}

The proposition is elementary, but it shows how the rescaling of~\cite{cartis2026note} is to be read. Passing from $\sfS^\intercal$ to $\sfS^\intercal/\sqrt d$ improves the evaluation bound by a factor of order $d/p$, from $\mcO(d^2\epsilon^{-2})$ to $\mcO(pd\epsilon^{-2})$ in the accounting of Table~\ref{tab:dim}, and the improvement is obtained with the trust-region, criticality and Hessian-bound parameters held fixed in subspace coordinates; by the proposition this is the same as running the unscaled algorithm with an ambient trust region smaller by $\sqrt d$, a criticality constant $\mu$ smaller by $d$ and a Hessian bound larger by $d$. Neither reading is wrong, but a table of complexity against sketch scale is a table over parametrizations, and it is meaningful based on which parameters are held fixed. We follow~\cite{cartis2026note} and Part~I and hold $(\Delta_0,\Delta_{\max},\mu,\kappa_s,\kH,\kedF)$ fixed in subspace coordinates; Section~\ref{sec:dimension} is written under that convention. What pins the scale in the bounds is the additive term coming from $\mu^{-1}$ and from the floor $\kH\ge1$, which do not rescale.

\section{Estimation error at a stopping-time sample size}\label{sec:accuracy}

This section contains the probabilistic content concerning the estimates. It has four parts: the sample size is an almost surely finite stopping time (Section~\ref{sec:finiteness}); a Bernstein inequality holds at that stopping time by optional stopping (Section~\ref{sec:stopped}); the pointwise route of~\cite{ha2025complexity} has two steps that read differently  from the stopping-time viewpoint (Section~\ref{sec:repairs}); and the bound is made uniform over the trust region, which is what admits reuse and what Borel--Cantelli then turns into an eventual statement (Section~\ref{sec:uniform}).

\subsection{The sample size is a finite stopping time}\label{sec:finiteness}

\begin{lemma}[freezing]\label{lem:freeze}
Conditionally on $\mcF^{\sfQ}_{k-1}$, the radius $\Delta_k$, the sketch $\sfQ_k$, the design points $\{\BFX_k^i\}_{i\le 2p_k}$ and the floor $n_k$ are fixed, and the replications $\xi_{k,1},\xi_{k,2},\dots$ are i.i.d.\ with the law of $\xi$. In particular, for each $i$, $\{F(\BFX_k^i,\xi_{k,j})\}_{j\ge1}$ are i.i.d.\ under the conditional law with mean $f(\BFX_k^i)$ and variance $\sigma_F^2(\BFX_k^i)\le v_F$.
\end{lemma}

\begin{proof}
The design points are $\mcF^{\sfQ}_{k-1}$-measurable, and $\sigma(\xi_{k,1},\xi_{k,2},\dots)$ is independent of $\mcF^{\sfQ}_{k-1}$ because the replications of iteration $k$ are fresh. The conditional law of a family independent of a $\sigma$-field is its unconditional law.
\end{proof}

\begin{proposition}[the sample size is a finite stopping time]\label{prop:Nkfinite}
Under Assumption~\ref{assum:variance}, $N_k$ is a stopping time with respect to $\{\mcG_{k,n}\}_{n\ge0}$, satisfies $N_k\ge n_k$ with $n_k$ as in~\eqref{eq:floor}, and $\mbP\bigl(\bigcap_{k}\{N_k<\infty\}\bigr)=1$.
\end{proposition}

\begin{proof}
For each $n$ the criterion in~\eqref{eq:stopping} is a function of $\Delta_k$, $\lambda_k$, the design points and $\xi_{k,1},\dots,\xi_{k,n}$, so $\{N_k\le n\}\in\mcG_{k,n}$. The floor is~\eqref{eq:floor}. For finiteness, fix $k$ and argue under the conditional law given $\mcF^{\sfQ}_{k-1}$ of Lemma~\ref{lem:freeze}: the ordinary strong law of large numbers gives $\sigmahat_F(\BFX_k^i,n)\to\sigma_F(\BFX_k^i)$ as $n\to\infty$ for each of the $2p_k+1$ indices, and a maximum over finitely many almost surely convergent sequences converges, so $$\sigma_0\vee\sigmahat^{\max}_F(k,n)\to\sigma_0\vee\max_i\sigma_F(\BFX_k^i)\le\sigma_0\vee\sqrt{v_F}<\infty.$$ The left side of the criterion tends to zero while the right side is a fixed positive number, so the criterion is met for all large $n$ and $N_k<\infty$ under the conditional law; integrating gives $\mbP(N_k<\infty)=1$, and a countable intersection of almost sure events is almost sure.
\end{proof}

Finiteness matters for two reasons beyond bookkeeping. The estimates $\Fbar(\BFX_k^i,N_k)$ presuppose a finite $N_k$, so without the proposition the algorithm might, with positive probability, sample forever and never reach the subproblem. And the optional stopping argument of the next subsection is an argument about $L_{N_k}$, which exists only when $N_k<\infty$; this is the point at which the present route uses finiteness explicitly, see Section~\ref{sec:repairs}. The conditional strong law of Majerek, Nowak and Zi\c{e}ba used in~\cite{dzahini2024stochastic} is not needed: the freezing lemma makes the replications i.i.d.\ under the conditional law, and the ordinary strong law applies. The rule~\eqref{eq:stopping} is, conditionally, a purely sequential rule of Chow--Robbins type~\cite{ChowRobbins1965}, and its finiteness and moments are classical~\cite{Mukhopadhyay1978,GhoshMukhopadhyay1981,GhoshMukhopadhyaySen1997,Shashaani2016Thesis}; only the finiteness is used here.

\subsection{Bernstein's inequality at the stopping time, by optional stopping}\label{sec:stopped}

\begin{lemma}[stopped Bernstein inequality]\label{lem:stopped}
Let $\{\mcG_j\}_{j\ge0}$ be a filtration and $\{Z_j\}_{j\ge1}$ be adapted to it ($Z_j$ is $\mcG_j$-measurable) with  $\mbE[Z_j\mid\mcG_{j-1}]=0$ and, for every integer $m\ge2$,
\begin{equation}
\mbE\bigl[|Z_j|^m\mid\mcG_{j-1}\bigr]\le\frac{m!}{2}B^{m-2}V^2\qquad\text{almost surely},\label{eq:momentCond}
\end{equation}
with $B,V>0$ $\mcG_0$-measurable. Let $A>0$ and a positive integer $\aleph$ be $\mcG_0$-measurable, and let $N$ be a stopping time for $\{\mcG_j\}$ with $\aleph\le N<\infty$ almost surely. Then, almost surely,
\begin{equation}
\mbP\Bigl(\frac1N\Bigl|\sum_{j=1}^N Z_j\Bigr|\ge A\ \Big|\ \mcG_0\Bigr)\le2\exp\Bigl(-\frac{\aleph A^2}{2(V^2+BA)}\Bigr).\label{eq:stopped}
\end{equation}
\end{lemma}

\begin{proof}
Set $S_n:=\sum_{j\le n}Z_j$ and $\Upsilon:=A/(V^2+BA)$, a $\mcG_0$-measurable random variable with $0<\Upsilon<B^{-1}$. For $z\in\mbR$, $$e^{\Upsilon z}\le1+\Upsilon z+\sum_{m\ge2}\Upsilon^m|z|^m/m!.$$ Taking conditional expectations and using~\eqref{eq:momentCond} and monotone convergence,
\[
\mbE\bigl[e^{\Upsilon Z_j}\mid\mcG_{j-1}\bigr]\le1+\frac{\Upsilon^2V^2}{2}\sum_{m\ge2}(\Upsilon B)^{m-2}=1+\frac{\Upsilon^2V^2}{2(1-\Upsilon B)}\le e^{\psi}\] where \[\psi:=\frac{\Upsilon^2V^2}{2(1-\Upsilon B)}=\frac{A^2}{2(V^2+BA)} .
\]
Hence $L_n:=\exp(\Upsilon S_n-n\psi)$ is a nonnegative supermartingale with $L_0=1$: $$\mbE[L_n\mid\mcG_{n-1}]=L_{n-1}e^{-\psi}\mbE[e^{\Upsilon Z_n}\mid\mcG_{n-1}]\le L_{n-1}.$$ For each $n$, $N\wedge n$ is a bounded stopping time, so the optional stopping theorem for supermartingales gives $$\mbE[L_{N\wedge n}\mid\mcG_0]\le L_0=1.$$ Since $N<\infty$ almost surely, $L_{N\wedge n}\to L_N$ almost surely, and Fatou's lemma gives $\mbE[L_N\mid\mcG_0]\le1$. On the event $\{S_N\ge NA\}$,
\[
L_N\ge\exp\bigl(N(\Upsilon A-\psi)\bigr)=\exp\Bigl(N\frac{A^2}{2(V^2+BA)}\Bigr)\ge\exp\Bigl(\aleph\frac{A^2}{2(V^2+BA)}\Bigr),
\]
using $N\ge\aleph$ and $\Upsilon A-\psi=A^2/(2(V^2+BA))>0$. The conditional Markov inequality then gives $$\mbP(S_N\ge NA\mid\mcG_0)\le e^{-\aleph A^2/(2(V^2+BA))}.$$ The same argument applied to $\{-Z_j\}$ bounds $\mbP(S_N\le-NA\mid\mcG_0)$, and the two together give~\eqref{eq:stopped}.
\end{proof}

\begin{remark}[what the stopping time is used for, and what it is not]\label{rem:stoppedVsVille}
Three statements about $S_n/n$ are available from the exponential supermartingale, at increasing strength and with the same constants. At a fixed $n\ge\aleph$: Bernstein's inequality. At the stopping time $N\ge\aleph$: Lemma~\ref{lem:stopped}, by optional stopping. Uniformly in $n\ge\aleph$: Ville's inequality~\cite{HowRamMcASek2020}, $$\mbP(\sup_{n\ge\aleph}|S_n|/n\ge A\mid\mcG_0)\le2e^{-\aleph A^2/(2(V^2+BA))},$$ which Part~I uses. The algorithm evaluates its estimates at exactly one sample size, $N_k$, and $N_k$ is a stopping time; Lemma~\ref{lem:stopped} is the statement that uses this property and no other. It requires $N<\infty$, which Proposition~\ref{prop:Nkfinite} supplies and which Ville's inequality does not need, and it claims nothing at sample sizes the algorithm never uses. We use it here for that reason and to keep this part independent of the maximal inequality; the choice has no effect on constants. What the lemma does \emph{not} give is a deviation scale adapted to $N$: the bound is at the scale $1/\sqrt{\aleph}$ of the floor, not $1/\sqrt N$. A scale $\sqrt{\log\log N/N}$ at the stopping time is available through the method of mixtures~\cite{Robbins1970,HowRamMcASek2020,delaPenaKlassLai2004}, at a $\log\log$ cost; it is not needed here because $N_k/n_k$ is bounded in probability by $(\sigma_0\vee\sqrt{v_F})^2/\sigma_0^2$ under the rule~\eqref{eq:stopping}, so the floor loses only a constant.
\end{remark}

\subsection{Two steps of the pointwise route, seen from the stopping time}\label{sec:repairs}
The estimation-error theorem of~\cite{ha2025complexity} asserts, for its zeroth-order common-random-number cases (B-0) and (C-0), that $\mbP(|\Ebar_k^{\text{s}}-\Ebar_k^0|>\kedF\Delta_k^2\ \text{i.o.})=0$. That is the conclusion needed here, and Corollary~\ref{cor:eventually} delivers it. The argument there reaches it through a pointwise inequality at a fixed sample size, summed over sample sizes. Two of its steps read differently once the sample size is treated as the stopping time, and we set them out. This is where the present analysis departs from its predecessor.

\begin{proposition}[the union bound over sample sizes]\label{prop:unionBound}
Let $Z_j$, $B$, $V$, $\aleph$ be as in Lemma~\ref{lem:stopped}, and let $A>0$. Bounding $\mbP(\sup_{n\ge\aleph}|S_n|/n\ge A\mid\mcG_0)$ by $\sum_{n\ge\aleph}\mbP(|S_n|/n\ge A\mid\mcG_0)$ and applying Bernstein's inequality at each $n$ gives
\begin{equation}
\sum_{n\ge\aleph}2e^{-nc}=\frac{2e^{-\aleph c}}{1-e^{-c}},\qquad c:=\frac{A^2}{2(V^2+BA)} .\label{eq:unionBound}
\end{equation}
In the application of~\cite{ha2025complexity}, $A=\kedF\Delta_k^2$, $V\asymp\Delta_k$ and $B\asymp\Delta_k$ (case (C-0)), so $c=c_k\asymp\Delta_k^2$, while $n_k\asymp\lambda_k\Delta_k^{-2}$. Hence $e^{-n_kc_k}=e^{-\Theta(\lambda_k)}$ is summable, but $$(1-e^{-c_k})^{-1}\ge c_k^{-1}\asymp\Delta_k^{-2}\to\infty$$ almost surely, and the product is not summable in general. Bounding the prefactor through the worst-case contraction $\Delta_k\ge\gamma_{dec}^k\Delta_0$, as in~\cite{ha2025complexity}, gives $(1-e^{-x^k})^{-1}\sim x^{-k}$ for $x<1$, which grows geometrically in $k$; a stopped inequality is the natural way to keep the bound summable.
\end{proposition}

\begin{proof}
Only the last sentence requires comment: $1-e^{-y}\le y$ for $y\ge0$, so $(1-e^{-c_k})^{-1}\ge c_k^{-1}$, and $c_k\to0$ because $\Delta_k\to0$ almost surely (Lemma~\ref{lem:squareSummable}).
\end{proof}


The point is not cosmetic: along the path, $\Delta_k^{-2}$ grows as fast as the sample size itself. Its origin is that a union bound treats every $n\ge n_k$ as a sample size the algorithm might use, whereas the algorithm uses exactly one, $N_k$, and that one is a stopping time. Lemma~\ref{lem:stopped} bounds the deviation at $N_k$ with the single exponential $2e^{-n_kc_k}$ and no prefactor, because optional stopping replaces the sum over $n$ by the evaluation of one supermartingale at one random time; Ville's inequality does the same with a supremum. Either serves the pointwise statements of~\cite{ha2025complexity} for every point that is fixed before sampling, namely the design points in all cases and the trial point whenever it is evaluated on replications drawn after the model is built.

\begin{proposition}[the trial point is not fixed before sampling]\label{prop:trialDependence}
Under common random numbers with reuse, $\BFX_k^{\text{s}}$ is a measurable function of $(\mcF^{\sfQ}_{k-1},\xi_{k,1},\dots,\xi_{k,N_k})$. Consequently the summands $D(\BFX_k^0,\BFX_k^{\text{s}};\xi_{k,j})$, $j\le N_k$, need not be adapted to $\{\mcG_{k,j}\}$, need not satisfy $\mbE[D(\BFX_k^0,\BFX_k^{\text{s}};\xi_{k,j})\mid\mcG_{k,j-1}]=0$, and the hypotheses of Lemma~\ref{lem:stopped}, of Ville's inequality and of Bernstein's inequality for martingales are not available.
\end{proposition}

\begin{proof}
$\BFX_k^{\text{s}}=\BFX_k^0+\sfQ_k\BFShat_k^*$ and $\BFShat_k^*$ is computed from $\Mhat_k$, which is a function of the estimates at the design points, which are functions of $\xi_{k,1},\dots,\xi_{k,N_k}$. For $j\le N_k$ the summand therefore depends on $\xi_{k,j+1},\dots,\xi_{k,N_k}$ through $\BFX_k^{\text{s}}$, and its conditional mean given $\mcG_{k,j-1}$ is $\mbE_{\xi}[D(\BFX_k^0,\BFy;\xi)]$ evaluated not at a fixed $\BFy$ but at a $\BFy$ correlated with $\xi_{k,j}$, which is not zero in general.
\end{proof}

In cases (B-0) and (C-0) of~\cite{ha2025complexity} the candidate shares the sample size $N_k$ of the design, and whether a pointwise inequality applies to it depends on which replications it is evaluated on. If the candidate and the incumbent are evaluated together on a stream drawn after the model is built, the pointwise inequality applies conditionally on the $\sigma$-field at the end of the model step, under which $\BFX_k^{\text{s}}$ is fixed, and the within-stream pairing keeps the variance proxy at the order that gives the rates of that paper. Cases (A-0), (A-1) and (C-1) use independent streams or a deterministic sample size, and case (B-1) needs only $\mcO(\Delta_k)$ accuracy along the segment, which a single net and a union bound supply uniformly (Remark~\ref{rem:singleNet}). Under full reuse, the scheme of Algorithm~\ref{alg:main}, no pointwise statement is available at the trial point, and what is needed instead is a bound holding simultaneously at every point of the trust region, which is indifferent to how the point was selected. The following remark places the three schemes side by side; the uniform bound is the content of the next subsection.

\begin{remark}[common random numbers are both the cause and the cure]\label{rem:no-crn}
It is worth being precise about what creates the need for a bound uniform in $\BFshat$, since the answer is neither ``common random numbers'' nor ``a random evaluation point'' by itself. Lemma~\ref{lem:stopped} requires that the threshold $A$ and the parameters $B,V$ be measurable with respect to the conditioning $\sigma$-field $\mcG_0$ and that the summands be a martingale difference sequence for $\{\mcG_j\}$. A random evaluation point is harmless if it is $\mcG_0$-measurable; what is fatal is a point that is a function of the very replications being averaged. Three sampling schemes illustrate this.

\emph{(i) No common random numbers.} Each design point and the candidate are evaluated on mutually independent streams. Let $\mathscr H_k$ be the $\sigma$-field generated by $\mcF^{\sfQ}_{k-1}$ and the replications used to build $\Mhat_k$. Then $\Mhat_k$, $\BFShat_k^*$ and $\BFX_k^{\text{s}}$ are $\mathscr H_k$-measurable and the candidate replications are independent of $\mathscr H_k$, so Lemma~\ref{lem:stopped} applies with $\mcG_0=\mathscr H_k$ and no uniformity in $\BFshat$ is needed. It is also unavailable: without shared randomness the field $\BFshat\mapsto\Dbar_k(\BFshat,n)$ of~\eqref{eq:Dbar} is not defined, and its natural substitute, the difference of two averages on independent streams, is as a function of $\BFshat$ a family of independent variables with no modulus of continuity. There is nothing to chain. The scheme pays for this in the variance: the difference of estimation errors between incumbent and candidate has standard deviation of order $\sigma_F/\sqrt n$ with no factor $\Delta_k$, so the tolerance $\kedF\Delta_k^2$ costs $n\asymp\Delta_k^{-4}$ in place of the $\Delta_k^{-2}$ of~\eqref{eq:floor}.

\emph{(ii) Common random numbers, with a separate stream for the candidate.} The design points share one stream; the incumbent and the candidate are then re-evaluated together on a second, independent stream. Pairing is preserved within each stream, so the variance proxy is again of order $\kQ^2\Delta_k^2$ and $n\asymp\Delta_k^{-2}$ suffices; and $\BFX_k^{\text{s}}$ is measurable with respect to the $\sigma$-field generated by the first stream, from which the second is independent, so Lemma~\ref{lem:stopped} applies as in (i). The cost is in the count, not the exponent: the incumbent is evaluated twice, $(2p_k+1)N_k+2N_k^{\text{s}}$ calls in place of $(2p_k+2)N_k$. Re-evaluating the incumbent is not optional, since pairing an estimate from the second stream with one from the first would break the coupling and return the variance proxy to order one. This is the scheme of ~\cite{ha2025complexity} and of the supermartingale analyses~\cite{chen2018stochastic,dzahini2024stochastic}, and Table~\ref{tab:schemes} shows why: it is the cheapest scheme for which a pointwise inequality suffices.

\emph{(iii) Full reuse, as in Algorithm~\ref{alg:main}.} One stream serves everything. Pairing is retained, the extra evaluations are avoided, and $\BFX_k^{\text{s}}$ is a function of the replications being averaged. There is no $\sigma$-field containing $\BFX_k^{\text{s}}$ from which the replications are independent, so Lemma~\ref{lem:stopped} cannot be applied at $\BFShat_k^*$, and a bound holding simultaneously over the trust region is needed; Lemma~\ref{lem:chained} supplies one, and the continuity that makes it possible is exactly the continuity supplied by the shared randomness.

Uniformity in $\BFshat$ is therefore required precisely when the evaluation point fails to be measurable with respect to a $\sigma$-field from which the replications are independent. Common random numbers within an iteration make the error field continuous, hence chainable; reuse across model and candidate destroys the measurability, hence makes chaining necessary. Scheme (i) has neither property, (ii) the first only, (iii) both. The saving of (iii) over (ii) is $2N_k^{\text{s}}-N_k$ calls, one evaluation of the incumbent when the two rules stop at the same size, a factor $(2p_k+3)/(2p_k+2)$ in the oracle count, and, more to the point of the analysis, the removal of a second sampling phase, a second stopping rule and a third $\sigma$-field. Whether that is worth the chaining argument is a matter of taste; that (iii) is what an implementation does is not.
\end{remark}

\begin{table}[h]\centering\small
\caption{Sampling schemes at iteration $k$. ``Uniformity'' refers to whether the accuracy statement must hold simultaneously over $\|\BFshat\|\le\Delta_k$.}\label{tab:schemes}
\begin{tabular}{llllll}
\toprule
scheme & streams & pairing & sample size & oracle calls & uniformity\\
\midrule
(i) no CRN & one per point & none & $\asymp\Delta_k^{-4}$ & $\sum_{i\le 2p_k}N^i_k+N^{\text{s}}_k$ & not needed\\
(ii) CRN, separate stream & two streams & within streams & $\asymp\Delta_k^{-2}$ & $(2p_k+1)N_k+2N^{\text{s}}_k$ & not needed\\
(iii) full reuse & one stream & throughout & $\asymp\Delta_k^{-2}$ & $(2p_k+2)N_k$ & needed\\
\bottomrule
\end{tabular}
\end{table}

\subsection{Uniform accuracy over the trust region, and the pathwise step}\label{sec:uniform}

The uniform statement is obtained by chaining: the unit ball is covered by nets at geometrically decreasing scales, the increment of $\Dbar_k(\cdot,n)$ along each edge is a paired increment by~\eqref{eq:telescope} and satisfies a Bernstein condition by Assumption~\ref{assum:bernstein}, and a union bound over the countably many edges controls them all at once. Part~I gives the argument in full (its Lemma~4.3). We only need it at the stopping time.

\begin{lemma}[chained Bernstein inequality at a stopping time]\label{lem:chained}
Let $\mathbb B_p:=\{\mathbf t\in\real^p:\|\mathbf t\|\le1\}$, $\{\mcG_j\}_{j\ge0}$ a filtration, $\Delta>0$ $\mcG_0$-measurable, and $\{Z_j(\mathbf t)\}_{j\ge1,\mathbf t\in\mathbb B_p}$ an adapted random field ($Z_j(\mathbf t)$ is $\mcG_j$-measurable) satisfying, almost surely: 
\begin{itemize}
    \item[(i)] $Z_j(\bm0)=0$ and $\mathbf t\mapsto Z_j(\mathbf t)$ is continuous on $\mathbb B_p$; 
    \item[(ii)] $\mbE[Z_j(\mathbf t)\mid\mcG_{j-1}]=0$ for every $\mathbf t\in\mathbb B_p$; 
    \item[(iii)] for all $\mathbf t,\mathbf t'\in\mathbb B_p$ with $\|\mathbf t-\mathbf t'\|\le1$ and integers $m\ge2$, $$\mbE[|Z_j(\mathbf t)-Z_j(\mathbf t')|^m\mid\mcG_{j-1}]\le\frac{m!}{2}(B\Delta\|\mathbf t-\mathbf t'\|)^{m-2}(V\Delta\|\mathbf t-\mathbf t'\|)^2$$ with $B,V>0$ $\mcG_0$-measurable. 
\end{itemize}
Let $\aleph$ be a positive integer $\mcG_0$-measurable random variable and $N$ a stopping time with $\aleph\le N<\infty$ almost surely. Then for every $\vartheta>0$, almost surely,
\begin{equation}
\mbP\Bigl(\sup_{\mathbf t\in\mathbb B_p}\frac1N\Bigl|\sum_{j=1}^N Z_j(\mathbf t)\Bigr|\ge15V\Delta\sqrt{\frac{p+\vartheta}{\aleph}}+72B\Delta\frac{p+\vartheta}{\aleph}\ \Big|\ \mcG_0\Bigr)\le e^{-\vartheta}.\label{eq:chained}
\end{equation}
\end{lemma}
\begin{proof}
Take from the proof of Part~I, Lemma~4.3, the nets $\mcN_m$ of $\mathbb B_p$ at scale $2^{-m}$ ($\mcN_0=\{\bm0\}$, $|\mcN_m|\le(3\cdot2^m)^p$), the nearest-point maps $\pi_m$, the sets $\mcE_m\subset\mcN_{m-1}\times\mcN_m$ of edges of length at most $\min\{3\cdot2^{-m},1\}$, which contain every pair $(\pi_{m-1}(\mathbf t),\pi_m(\mathbf t))$, and the levels $\vartheta_m:=\vartheta+4p(m+1)$, for which $4|\mcE_m|\le2^{-m}e^{\vartheta_m-\vartheta}$. 

Set $S_n(\mathbf t):=\sum_{j\le n}Z_j(\mathbf t)$ and $R_m:=3\cdot2^{-m}\Delta$. For an edge $(\mathbf u,\mathbf v)\in\mcE_m$ the increments $Z_j(\mathbf v)-Z_j(\mathbf u)$ are adapted, conditionally centered by~(ii) and satisfy~\eqref{eq:momentCond} with $BR_m$, $VR_m$ by~(iii), so Lemma~\ref{lem:stopped} with $A_m:=R_m\bigl(V\sqrt{2\vartheta_m/\aleph}+2B\vartheta_m/\aleph\bigr)$ gives $$\mbP(|S_N(\mathbf v)-S_N(\mathbf u)|\ge NA_m\mid\mcG_0)\le2e^{-\vartheta_m},$$ the exponent being at least $\vartheta_m$ by the computation in Part~I. A union bound over the countably many edges of all levels gives total conditional probability at most $\sum_{m\ge1}2|\mcE_m|e^{-\vartheta_m}\le e^{-\vartheta}$. Outside that event, for every $\mathbf t$ and $\phi\ge1$, telescoping along the edges $(\pi_{m-1}(\mathbf t),\pi_m(\mathbf t))$, $m\le\phi$, and using $S_N(\bm0)=0$ gives $|S_N(\pi_\phi(\mathbf t))|\le N\sum_{m\le\phi}A_m$; $S_N(\cdot)$ is continuous by~(i), $\pi_\phi(\mathbf t)\to\mathbf t$, and $$\sum_{m\ge1}A_m\le15V\Delta\sqrt{(p+\vartheta)/\aleph}+72B\Delta(p+\vartheta)/\aleph$$ by Part~I's summation. So here, only the time $N$ enters: no statement at any other sample size is used.
\end{proof}

\begin{remark}[why a single net does not suffice]\label{rem:singleNet}
One might try to avoid chaining by covering the ball with one net at the target scale and applying Lemma~\ref{lem:stopped} at each net point. Since the target accuracy is $\kedF\Delta_k^2$ and the process is Lipschitz in $\BFshat$ with modulus of order $\kQ$, the net must have mesh $\asymp\Delta_k$ in the coordinates where the accuracy is $\Delta_k^2$, that is, $\Delta_k^{-p}$ points; the union bound then costs $\log(\Delta_k^{-p})=p\log(1/\Delta_k)$ in the exponent, so the sample size would need $\lambda_k\gtrsim p\log(1/\Delta_k)\asymp pk$ for a radius contracting geometrically. Chaining pays $\sqrt{p}$ in a constant instead. The difficulty arises precisely because the accuracy threshold $\Delta_k^2$ is one power of $\Delta_k$ finer than the modulus $\Delta_k$ of the process; a threshold of order $\Delta_k$ is met by a net of $\mcO(1)^p$ points and a union bound. This is why only the zeroth-order regimes that both use common random numbers and require $\mcO(\Delta_k^2)$ accuracy of the estimate differences are affected: the first-order regime needs only $\mcO(\Delta_k)$ along the segment.
\end{remark}

We now apply Lemma~\ref{lem:chained} to the paired increment. Define
\begin{equation}
\rF:=\min\Bigl\{1,\ \Bigl(\frac{\sigma_0\kedF}{30\kappa_s\sigma_D\kQ}\Bigr)^2,\ \frac{\sigma_0^2\kedF}{144\kappa_s^2b_D\kQ\Delta_{\max}}\Bigr\},\qquad \cF:=\frac{\rF}{2},\label{eq:rF}
\end{equation}
and the accuracy event
\begin{equation}
\mcU_k:=\Bigl\{\sup_{\|\BFshat\|\le\Delta_k}|\Dbar_k(\BFshat,N_k)|\le\kedF\Delta_k^2\Bigr\}.\label{eq:Uk}
\end{equation}

\begin{theorem}[uniform weak accuracy; Part~I]\label{thm:uniform}
Let Assumptions~\ref{assum:variance}, \ref{assum:bernstein}, \ref{assum:lip-F} and~\ref{assum:well-aligned}(b) hold. If $\lambda_k\ge2p_k/\rF$, then almost surely
\begin{equation}
\mbP\bigl(\mcU_k^c\ \big|\ \mcF^{\sfQ}_{k-1}\bigr)\le e^{-\cF\lambda_k}.\label{eq:uniform}
\end{equation}
\end{theorem}

\begin{proof}
Apply Lemma~\ref{lem:chained} conditionally on $\mcF^{\sfQ}_{k-1}$ with $\mcG_j=\mcG_{k,j}$, $p=p_k$, $\Delta=\Delta_k$, $\aleph=n_k$, $N=N_k$ (a finite stopping time with $N_k\ge n_k$ by Proposition~\ref{prop:Nkfinite}) and the field $Z_j(\mathbf t):=D(\BFX_k^0,\BFX_k^0+\sfQ_k\Delta_k\mathbf t;\xi_{k,j})$ on $\mathbb B_{p_k}$, so that $\frac1n\sum_{j\le n}Z_j(\BFshat/\Delta_k)=\Dbar_k(\BFshat,n)$ for $\|\BFshat\|\le\Delta_k$. 

The field is adapted; (i) holds by~\eqref{eq:D}, Assumption~\ref{assum:lip-F} and the continuity of $f$; (ii) and (iii) hold with $B=b_D\kQ$ and $V=\sigma_D\kQ$, because by~\eqref{eq:telescope} $Z_j(\mathbf t)-Z_j(\mathbf t')$ is the paired increment at displacement $$\|\sfQ_k\Delta_k(\mathbf t-\mathbf t')\|\le\kQ\Delta_k\|\mathbf t-\mathbf t'\|\le\kQ\Delta_{\max}$$ for $\|\mathbf t-\mathbf t'\|\le1$, to which Assumption~\ref{assum:bernstein} applies conditionally on $\mcG_{k,j-1}$ by Lemma~\ref{lem:freeze}. Take $\vartheta=\cF\lambda_k$. Substituting $n_k\ge\sigma_0^2\lambda_k/(\kappa_s^2\Delta_k^2)$ and $\Delta_k\le\Delta_{\max}$ into~\eqref{eq:chained}, the threshold is at most
\[
\Delta_k^2\Bigl[\frac{15\kappa_s\sigma_D\kQ}{\sigma_0}\sqrt{\frac{p_k+\vartheta}{\lambda_k}}+\frac{72\kappa_s^2b_D\kQ\Delta_{\max}}{\sigma_0^2}\,\frac{p_k+\vartheta}{\lambda_k}\Bigr],
\]
and $(p_k+\vartheta)/\lambda_k\le\rF$ by $\lambda_k\ge2p_k/\rF$ and $\vartheta/\lambda_k=\rF/2$, so~\eqref{eq:rF} bounds each term by $\kedF\Delta_k^2/2$. Hence $\mcU_k^c$ is contained in the event of~\eqref{eq:chained}, whose conditional probability is at most $e^{-\cF\lambda_k}$. This is Part~I's Theorem~4.2 with the stopped lemma in place of the time-uniform one; only $\rF$, $\cF$ and $\lambda_k\ge2p_k/\rF$ are used below.
\end{proof}

The constants $\rF,\cF$ do not depend on $p_k$; the dependence on the subspace dimension sits entirely in the hypothesis $\lambda_k\ge2p_k/\rF$, which is the cost of uniformity over a $p_k$-dimensional ball (Section~\ref{sec:dimension}). Since $p_k=\mcO(\log k)$ while $\lambda_k/\log k\to\infty$, that hypothesis holds for all $k\ge k_\lambda$ with $k_\lambda$ deterministic; it is a condition on an index, not on $\lambda_0$, and this is one of the things $\lambda_k\to\infty$ buys.

\begin{corollary}[the pathwise step: accuracy holds eventually]\label{cor:eventually}
Under the hypotheses of Theorem~\ref{thm:uniform}, $\mbP(\mcU_k^c\ \text{i.o.})=0$. Equivalently, there is an almost surely finite random index $\Ku$ such that $\mcU_k$ holds for all $k\ge\Ku$, and on $\mcU_k$,
\begin{equation}
|\Ebar_k^i-\Ebar_k^0|\le\kedF\Delta_k^2\ \ (i=1,\dots,2p_k),\qquad |\Ebar_k^{\text{s}}-\Ebar_k^0|\le\kedF\Delta_k^2 .\label{eq:UkConsequences}
\end{equation}
\end{corollary}

\begin{proof}
For $k\ge k_\lambda$, taking expectations in~\eqref{eq:uniform} gives $$\mbP(\mcU_k^c)\le\exp(-\cF\lambda_0(\log(k+2))^{1+\epsilon_\lambda}),$$ which is summable in $k$ because $\cF\lambda_0(\log(k+2))^{1+\epsilon_\lambda}\ge2\log k$ for all large $k$, whatever $\lambda_0>0$. The first Borel--Cantelli lemma gives the claim; no independence among the $\mcU_k$ is used, and none is available, since $\mcU_k$ depends on $\BFX_k^0$ and $\Delta_k$ and hence on all previous replications. The consequences~\eqref{eq:UkConsequences} are the specializations of $\mcU_k$ to $\BFshat=\BFShat_k^i$ and $\BFshat=\BFShat_k^*$, both admissible since $\|\BFShat_k^i\|=\Delta_k$ and $\|\BFShat_k^*\|\le\Delta_k$.
\end{proof}

Corollary~\ref{cor:eventually} is the whole pathwise device in one statement: on almost every path there is a random, finite $\Ku(\omega)$ past which the estimates are accurate at every design point and at the trial point simultaneously, and the constant $\kedF$ can be any positive number (a smaller $\kedF$ only moves $\Ku$). Once alignment also holds eventually (Section~\ref{sec:alignment}), everything after this point is deterministic along the path. 
The corollary is also the form that~\cite[Theorem~3.2]{ha2025complexity} takes for cases (B-0) and (C-0) once the sample size is used as a stopping time: so no sum over $n$ or pointwise statement at the trial point is needed.

\begin{remark}[what the adaptive rule buys]\label{rem:ruleBuys}
Theorem~\ref{thm:uniform} uses the rule~\eqref{eq:stopping} only through the floor $N_k\ge n_k$, which follows from the guard $\sigma_0$ alone; no empirical standard deviation appears in any hypothesis. A rule calibrated instead to a global bound $\sigma_{F,\max}$ would meet the same accuracy guarantee while expending a factor $(\sigma_{F,\max}/\sigma_{k,+})^2$ more samples at every iteration, where $\sigma_{k,+}:=\sigma_0\vee\max_i\sigma_F(\BFX_k^i)$ is the noise level at the current design.
What the rule buys is efficiency: 
conditionally on $\mcF^{\sfQ}_{k-1}$ it is a Chow--Robbins rule and $N_k/(n_k\sigma_{k,+}^2/\sigma_0^2)\to1$ as $\lambda_k\to\infty$~\cite{ChowRobbins1965,GhoshMukhopadhyay1981}, so it attains the sample size an oracle knowing $\sigma_{k,+}$ would prescribe. This factor is invisible in an almost sure statement with a random constant and visible only in an expectation with deterministic constants, which is why it is displayed in Part~I and not here.
\end{remark}

\section{Model quality: the pathwise interface}\label{sec:model}

This section extracts from $\mcU_k$ the three deterministic inequalities that the trust-region argument consumes, and nothing else. Set
\begin{equation}
a:=\tfrac12\kLg\kQ^2+\kedF,\qquad \kg(p):=\sqrt{p}\,a,\qquad \kmf(p):=2\kg(p)+\kLg\kQ^2+\kH .\label{eq:constants}
\end{equation}

\begin{lemma}[reduced gradient]\label{lem:ffhat}
$\nabla\fhat_k(\BFshat)=\sfQ_k^\intercal\nabla f(\BFX_k^0+\sfQ_k\BFshat)$, and $\nabla\fhat_k$ is $\kLg\kQ^2$-Lipschitz on $\real^{p_k}$. Consequently, for $\|\BFshat\|\le\Delta_k$, $$|\fhat_k(\BFshat)-\fhat_k(\bm0)-\nabla\fhat_k(\bm0)^\intercal\BFshat|\le\frac12\kLg\kQ^2\|\BFshat\|^2.$$
\end{lemma}
\begin{proof}
The chain rule, then $$\|\sfQ_k^\intercal(\nabla f(\BFx_1)-\nabla f(\BFx_2))\|\le\kQ\kLg\|\sfQ_k(\BFshat_1-\BFshat_2)\|\le\kLg\kQ^2\|\BFshat_1-\BFshat_2\|,$$ then the integral form of Taylor's theorem.
\end{proof}

\begin{lemma}[model gradient at the center]\label{lem:center}
On $\mcU_k$, $\|\BFGhat_k-\nabla\fhat_k(\bm0)\|\le\kg(p_k)\Delta_k$.
\end{lemma}
\begin{proof}
Fix $i\in\{1,\dots,p_k\}$. By~\eqref{eq:GH}, $$[\BFGhat_k]_i=g_{k,i}+(\Ebar_k^i-\Ebar_k^{p_k+i})/(2\Delta_k)$$ with $$g_{k,i}:=(\fhat_k(\Delta_k\BFehat^i)-\fhat_k(-\Delta_k\BFehat^i))/(2\Delta_k)=\frac{1}{2\Delta_k}\int_{-\Delta_k}^{\Delta_k}\partial_i\fhat_k(t\BFehat^i)\,dt,$$ so by Lemma~\ref{lem:ffhat}
\[
|g_{k,i}-\partial_i\fhat_k(\bm0)|\le\frac{1}{2\Delta_k}\int_{-\Delta_k}^{\Delta_k}\kLg\kQ^2|t|\,dt=\tfrac12\kLg\kQ^2\Delta_k .
\]
On $\mcU_k$, by~\eqref{eq:UkConsequences}, $$|\Ebar_k^i-\Ebar_k^{p_k+i}|\le|\Ebar_k^i-\Ebar_k^0|+|\Ebar_k^{p_k+i}-\Ebar_k^0|\le2\kedF\Delta_k^2,$$ so $$|[\BFGhat_k]_i-\partial_i\fhat_k(\bm0)|\le(\frac12\kLg\kQ^2+\kedF)\Delta_k=a\Delta_k.$$ Summing the squares of the $p_k$ coordinates gives the claim.
\end{proof}

\begin{lemma}[trial-point difference]\label{lem:trial}
On $\mcU_k$, $|\Ebar_k^{\text{s}}-\Ebar_k^0|=|\Dbar_k(\BFShat_k^*,N_k)|\le\kedF\Delta_k^2$.
\end{lemma}
\begin{proof}
Immediate from~\eqref{eq:Uk}, since $\|\BFShat_k^*\|\le\Delta_k$.
\end{proof}

\begin{lemma}[model against function]\label{lem:mf}
On $\mcU_k$, for every $\|\BFshat\|\le\Delta_k$,
\begin{equation}
\bigl|\bigl(\Mhat_k(\BFshat)-\Mhat_k(\bm0)\bigr)-\bigl(\fhat_k(\BFshat)-\fhat_k(\bm0)\bigr)\bigr|\le\tfrac12\kmf(p_k)\Delta_k^2 .\label{eq:mf}
\end{equation}
\end{lemma}
\begin{proof} By Lemma~\ref{lem:ffhat}
$$\Mhat_k(\BFshat)-\Mhat_k(\bm0)=\BFGhat_k^\intercal\BFshat+\frac12\BFshat^\intercal\sfHhat_k\BFshat$$ and $$\fhat_k(\BFshat)-\fhat_k(\bm0)=\nabla\fhat_k(\bm0)^\intercal\BFshat+\varrho$$ with $|\varrho|\le\frac12\kLg\kQ^2\Delta_k^2$. The difference is at most $$\|\BFGhat_k-\nabla\fhat_k(\bm0)\|\Delta_k+\frac12\kH\Delta_k^2+\frac12\kLg\kQ^2\Delta_k^2,$$ and Lemma~\ref{lem:center} gives~\eqref{eq:mf}.
\end{proof}

\begin{lemma}[the clip is eventually inactive]\label{lem:hessian}
On $\mcU_k$ the unclipped second differences in~\eqref{eq:GH} are bounded in absolute value by $\kLg\kQ^2+2\kedF$; hence if $\kH\ge\kLg\kQ^2+2\kedF$ the clipping in~\eqref{eq:GH} has no effect for $k\ge\Ku$.
\end{lemma}
\begin{proof}
The second difference of $\fhat_k$ along $\BFehat^i$ divided by $\Delta_k^2$ is bounded by $\kLg\kQ^2$ by Lemma~\ref{lem:ffhat}, and the estimation part is $$(\Dbar_k(\Delta_k\BFehat^i,N_k)+\Dbar_k(-\Delta_k\BFehat^i,N_k))/\Delta_k^2,$$ bounded by $2\kedF$ on $\mcU_k$.
\end{proof}

Lemmas~\ref{lem:center}--\ref{lem:mf} are the entire interface between the probabilistic and the deterministic parts of the paper. No bound on $\nabla\Mhat_k(\BFshat)-\nabla\fhat_k(\BFshat)$ away from $\BFshat=\bm0$ is used, so the models are neither assumed nor shown to be fully linear; the supremum in Theorem~\ref{thm:uniform} is over a single scalar process and is there only because the point at which Lemma~\ref{lem:trial} is applied was selected by the replications being averaged.

\begin{remark}[the role of the quadratic term]\label{rem:quadratic}
$\sfHhat_k$ enters the analysis in exactly two places: through $\frac12\kH\Delta_k^2$ in Lemma~\ref{lem:mf}, and through $\|\sfHhat_k\|\vee1\le\kH$ in Assumption~\ref{assum:cauchy}. Neither uses what $\sfHhat_k$ approximates. The central-difference gradient~\eqref{eq:GH} is the same whether or not a diagonal Hessian is fitted, and in Lemma~\ref{lem:center} the pairs $\pm\Delta_k\BFehat^i$ cancel the quadratic term exactly, so $\kg$ is free of $\kH$. Every statement of Sections~\ref{sec:conv}--\ref{sec:complexity} therefore holds verbatim for the linear model $\Fbar_k^0+\BFGhat_k^\intercal\BFshat$, with $\kH$ replaced by $1$ in~\eqref{eq:constants} and~\eqref{eq:ksuc} and the subproblem solved exactly by $\BFShat_k^*=-\Delta_k\BFGhat_k/\|\BFGhat_k\|$, for which Assumption~\ref{assum:cauchy} holds with $\kappa_{\mathrm{fcd}}=1$. This is the stochastic form of the classical fact that first-order convergence of model-based trust-region methods needs gradient accuracy and a bounded Hessian, never Hessian accuracy~\cite[Ch.~10]{conn2009introduction}. Only $\mcO(\Delta_k)$ gradient accuracy is consumed, so forward differences on $p_k+1$ points would suffice for the theory; the symmetric design is kept because the second differences it yields cost nothing, having been paid for by the gradient. What the quadratic term buys is not visible in first-order rates: a step that uses curvature and need not be the Cauchy direction, a step that can be shorter than $\Delta_k$, and, under a Lipschitz Hessian, a central-difference bias of order $\Delta_k^2$ that a decoupled interpolation radius could exploit~\cite{grapiglia2023quadratic}.
\end{remark}

\begin{remark}[Hessians in random subspaces]\label{rem:subspaceHessian}
It is often observed that curvature is hard to maintain in random subspace methods: the reduced Hessian $\sfQ_k^\intercal\nabla^2f\,\sfQ_k$ lives in a subspace redrawn at every iteration, so curvature estimated at $k$ cannot be carried to $k+1$ by the least-change updates of full-space methods~\cite{powell2004newuoa}, and a full quadratic interpolation model in the subspace costs $(p+1)(p+2)/2$ points rather than $\mcO(p)$. Both observations are correct and neither affects the present model. The diagonal $\sfHhat_k$ is computed from the $2p_k+1$ points already used for the gradient, at no additional oracle cost and with no attempt to carry it across iterations; the analysis uses it only as a bounded perturbation, and Remark~\ref{rem:quadratic} shows the analysis asks nothing more of it. The difficulty would become a disadvantage only for a second-order guarantee or for curvature reuse, which we do not study here.
\end{remark}

\section{Convergence}\label{sec:conv}

Fix the sufficient-decrease constant of a successful iteration,
\begin{equation}
\ta:=\frac{\eta\kappa_{\mathrm{fcd}}}{2\mu}\min\Bigl\{1,\frac{1}{\mu\kH}\Bigr\},\qquad\text{and assume throughout}\qquad \kedF<\ta .\label{eq:ta}
\end{equation}
The condition is on $\kedF$, a design constant we are free to choose small, and it is the same as in~\cite{ha2025complexity}.

\begin{lemma}[radii are square summable]\label{lem:squareSummable}
Under the hypotheses of Corollary~\ref{cor:eventually}, Assumption~\ref{assum:cauchy} and~\eqref{eq:ta}
$\SDelta:=\sum_{k\ge0}\Delta_k^2<\infty$ almost surely; in particular $\Delta_k\to0$ almost surely.
\end{lemma}
\begin{proof}
Fix $\omega$ outside the null set of Corollary~\ref{cor:eventually}. If $k\ge\Ku$ is successful then $\rhohat_k\ge\eta$ and $\mu\|\BFGhat_k\|\ge\Delta_k$, so by Assumption~\ref{assum:cauchy} and $\|\sfHhat_k\|\vee1\le\kH$,
\[
\Fbar_k^0-\Fbar_k^{\text{s}}\ge\eta\bigl(\Mhat_k(\bm0)-\Mhat_k(\BFShat_k^*)\bigr)\ge\frac{\eta\kappa_{\mathrm{fcd}}}{2}\frac{\Delta_k}{\mu}\min\Bigl\{\Delta_k,\frac{\Delta_k}{\mu\kH}\Bigr\}=\ta\Delta_k^2 .
\]
Since $$\Fbar_k^0-\Fbar_k^{\text{s}}=f(\BFX_k^0)-f(\BFX_k^{\text{s}})-\Dbar_k(\BFShat_k^*,N_k),$$ Lemma~\ref{lem:trial} gives $$f(\BFX_k^0)-f(\BFX_k^{\text{s}})\ge(\ta-\kedF)\Delta_k^2$$ on every successful $k\ge\Ku$. Let $\mcK=\{k_1<k_2<\cdots\}$ be the successful iterations with $k_1\ge\Ku$. The iterate moves only on successes, so $\BFX_{k_i}^{\text{s}}=\BFX_{k_{i+1}}^0$ and, summing,
\begin{equation}
(\ta-\kedF)\sum_{k\in\mcK}\Delta_k^2\le f(\BFX_{k_1}^0)-f^*<\infty .\label{eq:sumSucc}
\end{equation}
Between consecutive successes the radius expands once by at most $\gamma_{inc}$ and then contracts geometrically, so $$\sum_{k=k_i+1}^{k_{i+1}}\Delta_k^2\le\gamma_{inc}^2\Delta_{k_i}^2/(1-\gamma_{dec}^2);$$ the iterations before $k_1$ are finitely many, and if $\mcK$ is finite the radii decrease geometrically after its last element. Hence $\SDelta<\infty$.
\end{proof}

\begin{lemma}[a small radius forces success]\label{lem:success}
Set
\begin{equation}
\ksuc(p):=\min\Bigl\{\mu,\ \frac1{\kH},\ \frac{\kappa_{\mathrm{fcd}}(1-\eta)}{2\bigl(\kmf(p)+2\kedF\bigr)}\Bigr\} .\label{eq:ksuc}
\end{equation}
On $\mcU_k$, if $\Delta_k\le\ksuc(p_k)\|\BFGhat_k\|$ then iteration $k$ is successful. Consequently $$\mbP\bigl(\{\Delta_k\le\ksuc(p_k)\|\BFGhat_k\|\}\cap\{k\text{ fails}\}\ \text{i.o.}\bigr)=0.$$
\end{lemma}
\begin{proof}
The criticality test $\mu\|\BFGhat_k\|\ge\Delta_k$ holds by the first term of~\eqref{eq:ksuc}. Write
\[
\rhohat_k-1=\frac{\bigl(f(\BFX_k^{\text{s}})-f(\BFX_k^0)\bigr)-\bigl(\Mhat_k(\BFShat_k^*)-\Mhat_k(\bm0)\bigr)+\Dbar_k(\BFShat_k^*,N_k)}{\Mhat_k(\BFShat_k^*)-\Mhat_k(\bm0)} .
\]
The numerator is bounded in absolute value by $(\frac12\kmf(p_k)+\kedF)\Delta_k^2$ by Lemmas~\ref{lem:mf} and~\ref{lem:trial}. For the denominator, $\Delta_k\le\|\BFGhat_k\|/\kH$ by the second term, so Assumption~\ref{assum:cauchy} gives $$\Mhat_k(\bm0)-\Mhat_k(\BFShat_k^*)\ge\frac12\kappa_{\mathrm{fcd}}\|\BFGhat_k\|\Delta_k\ge\frac{\kappa_{\mathrm{fcd}}}{2\ksuc(p_k)}\Delta_k^2.$$ Hence $$|\rhohat_k-1|\le\ksuc(p_k)(\kmf(p_k)+2\kedF)/\kappa_{\mathrm{fcd}}\le\frac12(1-\eta)$$ by the third term, so $\rhohat_k\ge\eta$. The last claim is Corollary~\ref{cor:eventually}.
\end{proof}

\subsection{Alignment holds eventually when the subspace dimension grows}\label{sec:alignment}

\begin{lemma}[alignment eventually]\label{lem:alignEventually}
If $\sum_k(1-\bQk)<\infty$ then $\mbP(\mcA_k^c\ \text{i.o.})=0$: there is an almost surely finite random index $\Ka$ with $\mcA_k$ for all $k\ge\Ka$.
\end{lemma}
\begin{proof} Observe that 
$$\mbP(\mcA_k^c)=\mbE[\mbP(\mcA_k^c\mid\mcF_{k-1})]\le1-\bQk$$ is  summable; so apply the first Borel--Cantelli lemma. The events are not independent, since $\nabla f(\BFX_k^0)$ depends on $\sfQ_0,\dots,\sfQ_{k-1}$, and the first lemma does not require it.
\end{proof}

For the ensembles of Section~\ref{sec:sketch}, $1-\bQk=e^{-c_{\sfQ}p_k}$ with $c_{\sfQ}=\alpha^2/4$ for orthonormal columns and a similar constant for Gaussian or hashing transforms~\cite{cartis2026note,dzahini2024stochastic}, so summability is achieved with
\begin{equation}
p_k=\min\Bigl\{d,\ \Bigl\lceil\frac{(1+\delta)\log(k+2)}{c_{\sfQ}}\Bigr\rceil\Bigr\}=\mcO(\log k),\qquad\delta>0,\label{eq:pk}
\end{equation}
which, with the capped bound of Proposition~\ref{prop:Haar},  gives $$1-\bQk\le(k+2)^{-(1+\delta)}$$ while $p_k<d$ and $0$ afterwards. A fixed $p$ gives a fixed $\bQ<1$, for which alignment need not hold eventually; 
what that costs is the subject of Section~\ref{sec:fixedp}.

\begin{lemma}[radius bounded below]\label{lem:radiusLB}
Set
\begin{equation}
\ka(p):=\frac{\gamma_{dec}\,\aQ\,\ksuc(p)}{1+\ksuc(p)\kg(p)},\label{eq:ka}
\end{equation}
which is nonincreasing in $p$. Under Assumptions~\ref{assum:smooth}--\ref{assum:well-aligned} with $\sum_k(1-\bQk)<\infty$ with $\{p_k\}$ nondecreasing and~\eqref{eq:ta}, almost surely there is a finite random index $K:=\Ku\vee\Ka$ such that, for every $\epsilon>0$ with $\Delta_K\ge\ka(p_K)\epsilon$ and every integer $\bar k\ge K$,
   \begin{equation}
   \|\nabla f(\BFX_k^0)\|\ge\epsilon\ \text{ for all }K\le k\le\bar k\quad\Longrightarrow\quad\Delta_k\ge\ka(p_k)\epsilon\ \text{ for all }K\le k\le\bar k.
   \end{equation}
\end{lemma}
\begin{proof}
Fix $\omega$ outside the null sets of Corollary~\ref{cor:eventually} and Lemma~\ref{lem:alignEventually}. The base case $k=K$ is the hypothesis on $\epsilon$. Suppose $\Delta_{k-1}\ge\ka(p_{k-1})\epsilon$ for some $K<k\le\bar k$ and, for contradiction, $\Delta_k<\ka(p_k)\epsilon$. Since $\ka(p_k)\le\ka(p_{k-1})$ we have $\Delta_k<\Delta_{k-1}$, so iteration $k-1$ failed and $\Delta_k=\gamma_{dec}\Delta_{k-1}$. As $k-1\ge\Ku$, Lemma~\ref{lem:success} gives $$\Delta_{k-1}>\ksuc(p_{k-1})\|\BFGhat_{k-1}\|.$$ As $k-1\ge\Ka$, $\mcA_{k-1}$ holds, so by Lemmas~\ref{lem:ffhat} and~\ref{lem:center}
\[
\|\BFGhat_{k-1}\|\ge\|\nabla\fhat_{k-1}(\bm0)\|-\kg(p_{k-1})\Delta_{k-1}\ge\aQ\epsilon-\kg(p_{k-1})\Delta_{k-1}.
\]
Combining, $$\Delta_{k-1}(1+\ksuc(p_{k-1})\kg(p_{k-1}))>\ksuc(p_{k-1})\aQ\epsilon\ \Rightarrow\ \Delta_{k-1}>\ka(p_{k-1})\epsilon/\gamma_{dec},$$ whence $$\Delta_k=\gamma_{dec}\Delta_{k-1}>\ka(p_{k-1})\epsilon\ge\ka(p_k)\epsilon,$$ which is a contradiction.
\end{proof}

The induction runs on every index $k-1\ge K$, and it needs alignment at every one of them; that is exactly why $\sum_k(1-\bQk)<\infty$ is the natural hypothesis of this route.

\begin{theorem}[strong consistency]\label{thm:consistency}
Under Assumptions~\ref{assum:smooth}--\ref{assum:well-aligned} with $\{p_k\}$ nondecreasing, $p_k\le d$, and $\sum_k(1-\bQk)<\infty$ (which~\eqref{eq:pk} achieves for the ensembles of Section~\ref{sec:sketch}), almost surely $$\lim_{k\to\infty}\|\nabla f(\BFX_k^0)\|=0.$$
\end{theorem}
\begin{proof}
Fix $\omega$ outside the relevant null sets.

Suppose $\|\nabla f(\BFX_k^0)\|\ge\epsilon$ for all $k\ge\bar k$ and some $\epsilon>0$. Enlarging $K$ to $K\vee\bar k$ and shrinking $\epsilon$ so that the base case of Lemma~\ref{lem:radiusLB} holds, $\Delta_k\ge\ka(p_k)\epsilon$ for all $k\ge K$. By~\eqref{eq:constants}, \eqref{eq:ksuc}, and \eqref{eq:ka}, since $\ka$ is nonincreasing and $p_k\le d$, $\ka(p_k)\ge\ka(d)>0$, so $\sum_{k\ge K}\Delta_k^2\ge\epsilon^2\sum_{k\ge K}\ka^2(d)=\infty$, contradicting Lemma~\ref{lem:squareSummable}. So $$\sum_{k\ge K}\Delta_k^2\ge\epsilon^2\sum_{k\ge K}\ka^2(p_k)=\infty,$$ contradicting Lemma~\ref{lem:squareSummable}. Hence $\liminf_k\|\nabla f(\BFX_k^0)\|=0$.

Now, suppose $\limsup_k\|\nabla f(\BFX_k^0)\|>2\epsilon>0$. Then there are infinitely many $t_j$ with $\|\nabla f(\BFX_{t_j}^0)\|>2\epsilon$ and, by the first part, for each $t_j$ a first $\ell_j>t_j$ with $\|\nabla f(\BFX_{\ell_j}^0)\|\le\epsilon$. On the window $t_j\le k<\ell_j$ we have $\|\nabla f(\BFX_k^0)\|>\epsilon$. 
Since $\kg(p_k)\le\kg(d)$ and $\Delta_k\to0$, $\kg(p_k)\Delta_k\to0$; 
hence for all large $j$ and all $k$ in the window, $$\|\BFGhat_k\|\ge\aQ\epsilon-\kg(p_k)\Delta_k\ge\frac12\aQ\epsilon$$ by Lemma~\ref{lem:center} and $\mcA_k$. On a successful $k$ in the window, Assumption~\ref{assum:cauchy}, Lemma~\ref{lem:trial} and $\rhohat_k\ge\eta$ give
\[
f(\BFX_k^0)-f(\BFX_k^{\text{s}})\ge\eta\bigl(\Mhat_k(\bm0)-\Mhat_k(\BFShat_k^*)\bigr)-\kedF\Delta_k^2\ge\frac{\eta\kappa_{\mathrm{fcd}}\aQ\epsilon}{4}\min\Bigl\{\Delta_k,\frac{\aQ\epsilon}{2\kH}\Bigr\}-\kedF\Delta_k^2\ge c_\epsilon\Delta_k
\]
with $c_\epsilon:=\eta\kappa_{\mathrm{fcd}}\aQ\epsilon/8$, once $$\Delta_k\le\min\{\aQ\epsilon/(2\kH),\,c_\epsilon/\kedF\},$$ which holds for large $j$ because $\Delta_k\to0$. The iterate moves only on successes and $\|\BFX_k^{\text{s}}-\BFX_k^0\|\le\kQ\Delta_k$, so
\[
\|\BFX_{\ell_j}^0-\BFX_{t_j}^0\|\le\kQ\sum_{\substack{t_j\le k<\ell_j\\k\ \text{successful}}}\Delta_k\le\frac{\kQ}{c_\epsilon}\bigl(f(\BFX_{t_j}^0)-f(\BFX_{\ell_j}^0)\bigr).
\]
By~\eqref{eq:sumSucc}, $f(\BFX_k^0)$ is nonincreasing along successes from $\Ku$ on and bounded below, hence convergent, so the right side tends to zero as $j\to\infty$. But $\|\nabla f(\BFX_{t_j}^0)-\nabla f(\BFX_{\ell_j}^0)\|>\epsilon$ and Assumption~\ref{assum:smooth} give $\|\BFX_{\ell_j}^0-\BFX_{t_j}^0\|>\epsilon/\kLg$, which is a contradiction.
\end{proof}

\section{Complexity}\label{sec:complexity}

Let $\Teps:=\inf\{k:\|\nabla f(\BFX_k^0)\|\le\epsilon\}$, finite almost surely by Theorem~\ref{thm:consistency}, and $\Weps:=\sum_{k<\Teps}(2p_k+2)N_k$.

\subsection{A first statement, with no machinery}\label{sec:firstStatement}

The following is the fixed-threshold bound that the pathwise route gives immediately.

\begin{theorem}[fixed-threshold iteration bound]\label{thm:fixedThreshold}
Under the hypotheses of Theorem~\ref{thm:consistency}, for every $\epsilon>0$ there is an almost surely finite random integer $\aleph_\epsilon$ such that
\begin{equation}
\mbP\bigl(\Teps\,\epsilon^2\le\SDelta+\aleph_\epsilon\,\epsilon^2\bigr)=1 .\label{eq:fixedThreshold}
\end{equation}
\end{theorem}
\begin{proof}
Every $k<\Teps$ has $\|\nabla f(\BFX_k^0)\|>\epsilon$. Split these indices into $$\mathscr L_\epsilon:=\{k<\Teps:\Delta_k<\epsilon\}\qquad \text{and} \qquad\mathscr G_\epsilon:=\{k<\Teps:\Delta_k\ge\epsilon\}.$$ By Theorem~\ref{thm:consistency} and Lemma~\ref{lem:squareSummable}, both $\|\nabla f(\BFX_k^0)\|\to0$ and $\Delta_k\to0$, so $\aleph_\epsilon:=|\{k:\|\nabla f(\BFX_k^0)\|>\epsilon,\ \Delta_k<\epsilon\}|<\infty$ almost surely and $|\mathscr L_\epsilon|\le\aleph_\epsilon$; and $\epsilon^2|\mathscr G_\epsilon|\le\sum_{k\in\mathscr G_\epsilon}\Delta_k^2\le\SDelta$. Adding gives~\eqref{eq:fixedThreshold}.
\end{proof}

The statement is honest about what the pathwise route gives for free: the constant $\aleph_\epsilon$ is finite for each $\epsilon$ but may depend on $\epsilon$, so~\eqref{eq:fixedThreshold} is not yet a rate. Making $\aleph_\epsilon$ independent of $\epsilon$ is exactly what the radius lower bound of Lemma~\ref{lem:radiusLB} does, at the price of $\ka^{-2}$ and of alignment at every iteration.

\subsection{Iteration and sample complexity}\label{sec:rates}

\begin{theorem}[iteration complexity]\label{thm:iter}
Under the hypotheses of Theorem~\ref{thm:consistency} there is an almost surely finite random variable $R>0$ such that, for every deterministic sequence $\epsilon_\ell\downarrow0$,
\begin{equation}
\mbP\Bigl(T_{\epsilon_\ell}\,\epsilon_\ell^2\,\ka^2(p_{T_{\epsilon_\ell}})>R\ \ \text{i.o.}\Bigr)=0 .\label{eq:iter}
\end{equation}
With $p_k$ as in~\eqref{eq:pk} this implies  $$\Teps=\mcO(\pi_a\epsilon^{-2}\log(1/\epsilon))=\Otilde(\pi_a\epsilon^{-2})$$ almost surely, where $$\pi_a:=\sup_p\ka^{-2}(p)/p$$ collects the dependence on the remaining constants.
\end{theorem}
\begin{proof}
Fix $\omega$ outside the relevant null sets, let $K$ be as in Lemma~\ref{lem:radiusLB}. If $T_\epsilon$ stays bounded as $\epsilon\downarrow0$ (an exact stationary point is reached), $T_\epsilon\epsilon^2\to0$ and there is nothing to prove; otherwise take $\epsilon$ small enough that $K<\Teps$ and the base case holds. For $K\le k<\Teps$, $\|\nabla f(\BFX_k^0)\|>\epsilon$, so $\Delta_k\ge\ka(p_k)\epsilon\ge\ka(p_{\Teps})\epsilon$ by monotonicity, and
\[
(\Teps-K)\,\ka^2(p_{\Teps})\epsilon^2\le\sum_{k=K}^{\Teps-1}\Delta_k^2\le\SDelta\ \Rightarrow\ \Teps\le\SDelta\ka^{-2}(p_{\Teps})\epsilon^{-2}+K.
\]
Since the first term tends to infinity as $\epsilon\downarrow0$, $K$ is absorbed and~\eqref{eq:iter} holds with $R:=\SDelta+K\ka^2(p_K)$. For the rate, $\Teps\epsilon^2\le R\pi_ap_{\Teps}$ and $p_{T_\epsilon}\le c\log(T_\epsilon+2)$; with $t:=\Teps\epsilon^2$ this means $t\le c\pi_a\log(t\epsilon^{-2})$, which forces $t=\mcO(\pi_a\log(1/\epsilon))$.
\end{proof}

\begin{lemma}[realized sample size]\label{lem:Nk}
Let Assumptions~\ref{assum:variance} and~\ref{assum:fourth} hold, and $\SDelta<\infty$ almost surely. Suppose further that $p_k/\lambda_k$ is bounded. There are a deterministic $n_W<\infty$ and an almost surely finite random index $\Kw$ such that $N_k\le n_W\lambda_k\Delta_k^{-2}$ for all $k\ge\Kw$.
\end{lemma}
\begin{proof}
Fix $\varepsilon'\in(0,1)$ and let $$n_k^\dagger:=\lceil(1+\varepsilon')\lambda_k\kappa_s^{-2}\Delta_k^{-2}(\sigma_0^2\vee c_4^{1/2})\rceil.$$ Since $\sigma_F^2(\BFx)\le c_4^{1/2}$ by Jensen, the event $\{N_k>n_k^\dagger\}$ forces $$\sigmahat_F^2(\BFX_k^i,n_k^\dagger)-\sigma_F^2(\BFX_k^i)>\varepsilon'\sigma_0^2$$ at one of the $2p_k+1$ design points. Conditionally on $\mcF^{\sfQ}_{k-1}$ (Lemma~\ref{lem:freeze}), a union bound and Chebyshev's inequality with $\mbE|\sigmahat_n^2-\sigma^2|^2\le c'/n$, valid under a fourth moment with $c'$ depending only on $c_4$, give
\[
\mbP\bigl(N_k>n_k^\dagger\ \big|\ \mcF^{\sfQ}_{k-1}\bigr)\le\frac{(2p_k+1)c'}{(\varepsilon'\sigma_0^2)^2n_k^\dagger}\le c''\,\frac{(2p_k+1)}{\lambda_k}\Delta_k^2 .
\]
The right side is $\mcF_{k-1}$-measurable so the same bound holds conditionally on $\mcF_{k-1}$ by the tower property. Moreover,  $(2p_k+1)/\lambda_k$ is bounded because $p_k=\mcO(\log k)$, so its sum over $k$ is at most a constant times $\SDelta<\infty$. L\'evy's extension of the Borel--Cantelli lemma for $\{\mcF_{k-1}\}$ gives $N_k\le n_k^\dagger$ for all $k\ge\Kw$ with $\Kw$ almost surely finite, and $n_k^\dagger\le n_W\lambda_k\Delta_k^{-2}$ after absorbing the ceiling with $\lambda_k\Delta_k^{-2}\ge\lambda_0\Delta_{\max}^{-2}$.
\end{proof}

\begin{theorem}[sample complexity]\label{thm:work}
Under the hypotheses of Theorem~\ref{thm:consistency} and Assumption~\ref{assum:fourth}, there is an almost surely finite random variable $R_W>0$ such that, for every deterministic $\epsilon_\ell\downarrow0$,
\begin{equation}
\mbP\Bigl(W_{\epsilon_\ell}>R_W\,(2p_{T_{\epsilon_\ell}}+2)\,\lambda_{T_{\epsilon_\ell}}\,\ka^{-4}(p_{T_{\epsilon_\ell}})\,\epsilon_\ell^{-4}\ \ \text{i.o.}\Bigr)=0 ,\label{eq:work}
\end{equation}
and hence $\Weps=\Otilde(\pi_a^2\epsilon^{-4})$ almost surely.
\end{theorem}
\begin{proof}
Fix $\omega$ and let $K_w:=K\vee\Kw$. If $T_\epsilon$ stays bounded as $\epsilon\downarrow0$ (an exact stationary point is reached), $T_\epsilon\epsilon^2\to0$ and there is nothing to prove; otherwise take $\epsilon$ small enough that $K<T_\epsilon$. Split $$\Weps=\sum_{k<K_w}(2p_k+2)N_k+\sum_{k=K_w}^{\Teps-1}(2p_k+2)N_k,$$ the first sum a finite random number $Q_w$. For $K_w\le k<\Teps$, Lemma~\ref{lem:Nk} and Lemma~\ref{lem:radiusLB} give $$N_k\le n_W\lambda_k\ka^{-2}(p_k)\epsilon^{-2},$$ and $p_k$, $\lambda_k$, $\ka^{-2}(p_k)$ are nondecreasing, so the second sum is at most $$(2p_{\Teps}+2)n_W\lambda_{\Teps}\ka^{-2}(p_{\Teps})\epsilon^{-2}\Teps.$$ Substituting $\Teps\le\SDelta\ka^{-2}(p_{\Teps})\epsilon^{-2}+K$ from the proof of Theorem~\ref{thm:iter} gives~\eqref{eq:work} with $R_W:=2n_W\SDelta+Q_w$, the additive terms being absorbed as $\epsilon\downarrow0$. For the rate, $p_{\Teps}$, $\lambda_{\Teps}$ and $\ka^{-2}(p_{\Teps})/\pi_a$ are polylogarithmic in $1/\epsilon$ by Theorem~\ref{thm:iter}.
\end{proof}

\begin{remark}[on the form of the statements]\label{rem:asScope}
Both results are almost sure, with a random constant valid for $\epsilon$ below a random threshold, and they are stated in the "infinitely often along a deterministic sequence" form because $\epsilon_0(\omega)$ is random. This is the natural currency of the route, and we do not present the absence of expected bounds as a deficiency: the accuracy guarantee itself is an eventual statement, and everything downstream is deterministic along a fixed path. $R$ and $R_W$ contain the burn-in indices. Proposition~\ref{prop:burninTail} gives $\Ku$ moments of every order and $\Ka$ moments of order below $\delta$; since $\SDelta\le a_1+a_2\Ku$ under the hypotheses of Theorem~\ref{thm:finiteTime}, $R$ is integrable when $\delta>1$, while nothing is claimed for $R_W$, as no tail bound for $\Kw$ is given. What prevents a bound in expectation is not the constant but the random threshold $\epsilon_0(\omega)$: the threshold-free form~\eqref{eq:noThreshold} carries the factor $\gamma_{dec}^{-2K}/\Delta_0^2$, which is not integrable. Part~I obtains expected bounds, for fixed $p$, from its per-iteration probability bound. 
\end{remark}

\subsection{The burn-in is not arbitrarily late: tails and a finite-time bound}\label{sec:tails}

A natural objection to almost sure statements is that they are asymptotic: the good events hold from some $K(\omega)$ on, and nothing seems to be said about how late $K(\omega)$ can be, whereas an expected bound controls the tail of $\Teps$ through Markov's inequality. The objection does not survive a look at how $K(\omega)$ was constructed. It is the last index at which an event with an explicitly summable failure probability fails, so its tail is bounded by the tail of the series one sums in Borel--Cantelli, which is deterministic and available.

\begin{proposition}[tail of the burn-in]\label{prop:burninTail}
Under the hypotheses of Theorem~\ref{thm:consistency}, for every $k\ge k_\lambda$ and $\Ku:=1+\sup\{k:\mcU_k^c\}$ and $\Ka:=1+\sup\{k:\mcA_k^c\}$,
\begin{equation}
\mbP(\Ku>k)\le\sum_{j\ge k}e^{-\cF\lambda_j},\qquad \mbP(\Ka>k)\le\sum_{j\ge k}(1-\beta_{\sfQ,j}),\label{eq:burninTail}
\end{equation} 
and $\mbP(\Ka>k)=0$ for $k\ge k_d:=\min\{k:p_k=d\}$ for the ensemble of Proposition~\ref{prop:Haar} with~\eqref{eq:pk}.  
With $\lambda_j=\lambda_0(\log(j+2))^{1+\epsilon_\lambda}$ the first series is at most $2k\,e^{-\cF\lambda_k}$ for all large $k$, which decays faster than any power of $k$, so $\mbE[\Ku^m]<\infty$ for every $m$; with $p_k$ as in~\eqref{eq:pk}, the second is at most $k^{-\delta}/\delta$, so $\mbE[\Ka]<\infty$ when $\delta>1$.
\end{proposition}
\begin{proof} By Theorem~\ref{thm:uniform}
$$\{\Ku>k\}=\{\exists j\ge k:\mcU_j^c\}\subseteq\bigcup_{j\ge k}\mcU_j^c$$ and $$\mbP(\mcU_j^c)\le e^{-\cF\lambda_j}$$ for $j\ge k_\lambda$; likewise for $\Ka$ with Lemma~\ref{lem:alignEventually}. Writing $e^{-\cF\lambda_j}=(j+2)^{-a_j}$ with $$a_j:=\cF\lambda_0(\log(j+2))^{\epsilon_\lambda}\to\infty,$$ for $k$ large enough that $a_k\ge2$ one has $$\sum_{j\ge k}(j+2)^{-a_j}\le\sum_{j\ge k}(j+2)^{-a_k}\le k^{-a_k}(1+k/(a_k-1))\le2k\,e^{-\cF\lambda_k}.$$ The moments follow from $\mbE[K^m]=\sum_k\mbP(K^m>k).$ For $\Ka$, $$\sum_{j\ge k}(j+2)^{-(1+\delta)}\le k^{-\delta}/\delta.$$
\end{proof}

The proposition converts the almost sure statements into finite-time ones. We record the iteration bound. 

\begin{theorem}[finite-time, high-probability iteration bound]\label{thm:finiteTime}
Let the hypotheses of Theorem~\ref{thm:consistency} hold and let in addition $f$ be $\kappa_{\mathrm Lf}$-Lipschitz. Set $c_\Delta:=\gamma_{inc}^2/(1-\gamma_{dec}^2)$ and
\[
a_1:=c_\Delta\Bigl(\frac{\Delta_0^2}{\gamma_{dec}^2}+\frac{f(\BFx_0)-f^*}{\ta-\kedF}\Bigr),\qquad a_2:=c_\Delta\Bigl(\Delta_{\max}^2+\frac{\kappa_{\mathrm Lf}\kQ\Delta_{\max}}{\ta-\kedF}\Bigr),
\]
both deterministic. Then, almost surely, for every $\epsilon>0$ and with $K=\Ku\vee\Ka$,
\begin{equation}
\Teps\le K+\bigl(a_1+a_2\Ku\bigr)\Bigl(\frac{\gamma_{dec}^{-2K}}{\Delta_0^2}+\frac{1}{\ka^2(p_{\Teps})\epsilon^2}\Bigr),\label{eq:noThreshold}
\end{equation}
and consequently, for every $\epsilon>0$ and every $k_0\ge k_\lambda$,
\begin{equation}
\mbP\Bigl(\Teps>k_0+(a_1+a_2k_0)\Bigl(\frac{\gamma_{dec}^{-2k_0}}{\Delta_0^2}+\frac{1}{\ka^2(d)\epsilon^2}\Bigr)\Bigr)\le\sum_{j\ge k_0}\Bigl(e^{-\cF\lambda_j}+1-\beta_{\sfQ,j}\Bigr).\label{eq:finiteTime}
\end{equation}
\end{theorem}
\begin{proof}
Fix $\omega$ outside the relevant null sets. First, $\SDelta\le a_1+a_2\Ku$: in the proof of Lemma~\ref{lem:squareSummable}, $$\SDelta\le c_\Delta(\Delta_0^2/\gamma_{dec}^2+\sum_{k\in\mcK}\Delta_k^2);$$ the successes before $\Ku$ number at most $\Ku$ and each contributes at most $\Delta_{\max}^2$; the successes from $\Ku$ on contribute at most $(f(\BFX_{k_1}^0)-f^*)/(\ta-\kedF)$ by~\eqref{eq:sumSucc}; and $$f(\BFX_{k_1}^0)\le f(\BFx_0)+\Ku\kappa_{\mathrm Lf}\kQ\Delta_{\max}$$ since each of the at most $\Ku$ moves before $\Ku$ has length at most $\kQ\Delta_{\max}$. Second, the induction of Lemma~\ref{lem:radiusLB} with the base case replaced by $\min\{\Delta_K,\ka(p_K)\epsilon\}$ gives $$\Delta_k\ge\min\{\Delta_K,\ka(p_k)\epsilon\}$$ for all $K\le k<\Teps$ and every $\epsilon>0$, with no threshold: if $$\Delta_k<\min\{\Delta_K,\ka(p_k)\epsilon\}\le\Delta_{k-1}$$ then $k-1$ failed and the contradiction is reached exactly as there. Since $\Delta_K\ge\gamma_{dec}^K\Delta_0$, summing $\Delta_k^2$ over $K\le k<\Teps$ gives $$(\Teps-K)\min\{\gamma_{dec}^{2K}\Delta_0^2,\ka^2(p_{\Teps})\epsilon^2\}\le\SDelta,$$ which is~\eqref{eq:noThreshold} after $1/\min\{a,b\}\le1/a+1/b$. On $\{K\le k_0\}$ the right side of~\eqref{eq:noThreshold} is at most the threshold in~\eqref{eq:finiteTime}, so the event there is contained in $\{K>k_0\}$, and Proposition~\ref{prop:burninTail} bounds its probability.
\end{proof}

We note three remarks here. First, \eqref{eq:finiteTime} is not asymptotic in any variable: for a given $\epsilon$ and a given $k_0$ it is a numerical statement about the run, with a deterministic right-hand side that can be evaluated from the algorithm's parameters and the constants $\cF$, $c_{\sfQ}$; the almost sure Theorem~\ref{thm:iter} is its limit as $k_0\to\infty$ along $\epsilon_\ell\downarrow0$. Second, the term $\gamma_{dec}^{-2k_0}/\Delta_0^2$ is the price of not knowing the radius at the end of the burn-in, and it is what prevents integrability of the constant in~\eqref{eq:noThreshold}; in the regime of interest, $\epsilon\ll\Delta_0\gamma_{dec}^{k_0}\ka$, it is dominated by the $\epsilon^{-2}$ term, and it disappears from the almost sure statement because the random threshold $\epsilon_0(\omega)$ absorbs it. Third, the tail in~\eqref{eq:finiteTime} is superpolynomial in $k_0$ on the estimation side and polynomial, of order $k_0^{-\delta}$, on the alignment side; the latter is set by the growth constant in~\eqref{eq:pk} and can be made as fast as one likes at a constant-factor cost in $p_k$. A Markov bound from an expected iteration bound, by contrast, gives a tail of order $t^{-1}$ in the deviation $t$; the expectation route gets sharper tails only through the high-probability arguments it also contains.
\section{Fixed subspace dimension: what the pathwise route can and cannot do}\label{sec:fixedp}

Sections~\ref{sec:conv}--\ref{sec:complexity} required $p_k=\mcO(\log k)$. This section keeps $p_k=p$ fixed, so that the misalignment probability $\dQ=1-\bQ$ is a fixed positive number, and asks what survives. The answer is quantitative and, we think, the most instructive point of the paper for a reader used to statements in expectation: the liminf and the complexity bounds survive with a worse exponent, and the loss is exactly the difference between the mean of a geometric random variable and the running maximum of a sequence of them. Part~I obtains $\epsilon^{-2}$ at fixed $p$ through its per-iteration probability bound; what follows is what the eventual-accuracy argument obtains without it.

Throughout the section $\kg,\kmf,\ksuc,\ka$ are the constants of~\eqref{eq:constants}, \eqref{eq:ksuc}, \eqref{eq:ka} at the fixed $p$, and we write
\begin{equation}
b_\epsilon:=\frac{\aQ\ksuc\,\epsilon}{1+\ksuc\kg}=\frac{\ka\epsilon}{\gamma_{dec}},\qquad c_\gamma:=\frac{\log\gamma_{inc}}{\log(1/\gamma_{dec})},\qquad \rho_{\sfQ}:=\frac{e(1+c_\gamma)}{c_\gamma}\,\dQ .\label{eq:fixedpConstants}
\end{equation}
The threshold $b_\epsilon$ is the radius below which alignment and accuracy force success when $\|\nabla f(\BFX_k^0)\|\ge\epsilon$: on $\mcU_k\cap\mcA_k$, Lemma~\ref{lem:center} gives $$\|\BFGhat_k\|\ge\aQ\epsilon-\kg\Delta_k,$$ and $\Delta_k\le b_\epsilon$ is equivalent to $$\Delta_k\le\ksuc(\aQ\epsilon-\kg\Delta_k),$$ so Lemma~\ref{lem:success} applies. The constant $c_\gamma$ converts expansions into units of contraction; $c_\gamma=1$ when $\gamma_{inc}=1/\gamma_{dec}$.

\begin{lemma}[misaligned iterations in a window]\label{lem:chernoff}
Let $I_i:=\ind\{\mcA_i^c\}$. Under Assumption~\ref{assum:well-aligned}(a) with $\sup_k(1-\bQk)\le\dQ$, for all integers $j\ge-1$ (with $\mcF_{-1}$ trivial), $j'\ge1$ and all real $t\ge j'\dQ$, almost surely
\begin{equation}
\mbP\Bigl(\sum_{i=j+1}^{j+j'}I_i\ge t\ \Big|\ \mcF_j\Bigr)\le\Bigl(\frac{ej'\dQ}{t}\Bigr)^{t}.\label{eq:chernoff}
\end{equation}
\end{lemma}
\begin{proof}
For $\theta>0$, $$\mbE[e^{\theta I_i}\mid\mcF_{i-1}]=1+\mbP(\mcA_i^c\mid\mcF_{i-1})(e^\theta-1)\le1+\dQ(e^\theta-1)\le\exp(\dQ(e^\theta-1)),$$ so by iterated conditioning $$\mbE[e^{\theta\sum_{i=j+1}^{j+j'}I_i}\mid\mcF_j]\le\exp(j'\dQ(e^\theta-1)).$$ Markov's inequality gives $$\mbP\left(\sum I_i\ge t\mid\mcF_j\right)\le\exp(-\theta t+j'\dQ(e^\theta-1)),$$ and the choice $e^\theta=t/(j'\dQ)\ge1$ gives $$\exp(-t\log(t/(j'\dQ))+t-j'\dQ)\le(ej'\dQ/t)^t.$$
\end{proof}

\begin{lemma}[depth of an excursion below $b_\epsilon$]\label{lem:depth}
Fix $\omega$ outside the null set of Corollary~\ref{cor:eventually}, integers $\bar k\ge\Ku$ and $k>\bar k$, an integer $r\ge0$, and $\epsilon>0$ with $\gamma_{inc}b_\epsilon\le\Delta_{\max}$. Suppose that $\|\nabla f(\BFX_i^0)\|\ge\epsilon$ for all $\bar k\le i<k$ and that $\Delta_k<\min\{\gamma_{dec}b_\epsilon,\Delta_{\bar k}\}\,\gamma_{dec}^{\,r}$. Then there is $j'\in\{1,\dots,k-\bar k\}$ such that the window $\{k-j',\dots,k-1\}$ contains more than
\begin{equation}
t_{j'}(r):=\frac{r+c_\gamma j'}{1+c_\gamma}\label{eq:tm}
\end{equation}
misaligned iterations.
\end{lemma}
\begin{proof}
Let $j_0:=\max\{i\in[\bar k,k):\Delta_i\ge b_\epsilon\}$ if this set is nonempty and $j_0:=\bar k-1$ otherwise, and let $j':=k-j_0-1$, so that the window is $\{j_0+1,\dots,k-1\}\subset[\bar k,k)$. In the first case $\Delta_{j_0+1}\ge\gamma_{dec}\Delta_{j_0}\ge\gamma_{dec}b_\epsilon$, since the radius decreases by at most the factor $\gamma_{dec}$ in one iteration; in the second $\Delta_{j_0+1}=\Delta_{\bar k}$. Either way
\[
\Delta_{j_0+1}\ge\min\{\gamma_{dec}b_\epsilon,\Delta_{\bar k}\}>\Delta_k,
\]
so $j_0+1<k$, that is, $j'\ge1$. For every $i$ in the window, $\Delta_i<b_\epsilon\le\Delta_{\max}/\gamma_{inc}$ by the choice of $j_0$, so a success multiplies the radius by exactly $\gamma_{inc}$ and a failure by $\gamma_{dec}$. With $J$ the number of failures in the window, $\Delta_k=\Delta_{j_0+1}\gamma_{dec}^{J}\gamma_{inc}^{\,j'-J}$, and the hypothesis gives $\gamma_{dec}^{J}\gamma_{inc}^{\,j'-J}<\gamma_{dec}^{\,r}$. Taking logarithms and dividing by $\log(1/\gamma_{dec})$,
\[
J-c_\gamma(j'-J)>r,\qquad\text{that is,}\qquad J>\frac{r+c_\gamma j'}{1+c_\gamma}=t_{j'}(r).
\]
Finally, every failure in the window is a misaligned iteration: for such $i$ we have $i\ge\bar k\ge\Ku$, so $\mcU_i$ holds, $\Delta_i<b_\epsilon$ and $\|\nabla f(\BFX_i^0)\|\ge\epsilon$, and by the discussion after~\eqref{eq:fixedpConstants} the iteration would be successful if $\mcA_i$ held.
\end{proof}

\begin{theorem}[fixed subspace dimension]\label{thm:fixedp}
Let the hypotheses of Corollary~\ref{cor:eventually}, Assumptions~\ref{assum:smooth}, \ref{assum:cauchy} and~\ref{assum:well-aligned} and~\eqref{eq:ta} hold with $p_k=p$ fixed and $\sup_k(1-\bQk)\le\dQ$, where $0<\dQ<c_\gamma/(e(1+c_\gamma))$, so that $\rho_{\sfQ}<1$. Fix $\delta'>0$ and set
\begin{equation}
\varrho:=\frac{(1+c_\gamma)(1+\delta')\log(1/\gamma_{dec})}{\log(1/\rho_{\sfQ})} .\label{eq:aExponent}
\end{equation}
Then almost surely there is a finite random index $K_2>\Ku$ such that, for every $\bar k\ge\Ku$, every $k\ge K_2$ with $k>\bar k$, and every $\epsilon>0$ with $\gamma_{inc}b_\epsilon\le\Delta_{\max}$,
\begin{equation}
\|\nabla f(\BFX_i^0)\|\ge\epsilon\ \text{ for all }\bar k\le i<k\quad\Longrightarrow\quad\Delta_k\ge\min\{\gamma_{dec}b_\epsilon,\Delta_{\bar k}\}\,\gamma_{dec}\,k^{-\varrho}.\label{eq:fixedpRadius}
\end{equation}
Consequently:
\begin{enumerate}
\item[(i)] if $\varrho\le\frac12$ then $\liminf_k\|\nabla f(\BFX_k^0)\|=0$ almost surely;
\item[(ii)] if $\varrho<\frac12$ then there is an almost surely finite $\bar{R}>0$ with $\mbP\bigl(T_{\epsilon_\ell}>\bar{R}\epsilon_\ell^{-2/(1-2\varrho)}\ \text{i.o.}\bigr)=0$ along every deterministic $\epsilon_\ell\downarrow0$, and under Assumption~\ref{assum:fourth}, $\Weps=\Otilde(\epsilon^{-4/(1-2\varrho)})$ almost surely.
\end{enumerate}
\end{theorem}
\begin{proof}
\emph{The radius bound.} For $k\ge1$ let $$r_k:=\lceil(1+c_\gamma)(1+\delta')\log k/\log(1/\rho_{\sfQ})\rceil,$$ so that $\gamma_{dec}^{\,r_k}\ge\gamma_{dec}\,k^{-\varrho}$, and let $\mcE_k$ be the event that for some $j'\in\{1,\dots,k\}$ the window $\{k-j',\dots,k-1\}$ contains at least $t_{j'}(r_k)$ misaligned iterations. Fix $j'$. Since $$t_{j'}(r_k)\ge c_\gamma j'/(1+c_\gamma)>j'\dQ,$$ Lemma~\ref{lem:chernoff} applies with $j=k-j'-1$ ($\mcF_{-1}$ being the trivial $\sigma$-field) and $t=t_{j'}(r_k)$, and since $$ej'\dQ/t_{j'}(r_k)\le e(1+c_\gamma)\dQ/c_\gamma=\rho_{\sfQ}<1,$$ the probability that the window contains at least $t_{j'}(r_k)$ misaligned iterations is at most $\rho_{\sfQ}^{\,t_{j'}(r_k)}$. Summing over $j'$,
\[
\mbP(\mcE_k)\le\rho_{\sfQ}^{\,r_k/(1+c_\gamma)}\sum_{j'\ge1}\rho_{\sfQ}^{\,c_\gamma j'/(1+c_\gamma)}\le\frac{k^{-(1+\delta')}}{1-\rho_{\sfQ}^{\,c_\gamma/(1+c_\gamma)}},
\]
which is summable. By the first Borel--Cantelli lemma there is an almost surely finite $K_2'$ with $\mcE_k^c$ for all $k\ge K_2'$; set $K_2:=K_2'\vee(\Ku+1)$. Now fix $\omega$ outside the null sets, $\bar k\ge\Ku$, $k\ge K_2$ with $k>\bar k$ and $\epsilon$ as in the statement, and suppose the hypothesis of~\eqref{eq:fixedpRadius} holds while $$\Delta_k<\min\{\gamma_{dec}b_\epsilon,\Delta_{\bar k}\}\gamma_{dec}k^{-\varrho}\le\min\{\gamma_{dec}b_\epsilon,\Delta_{\bar k}\}\gamma_{dec}^{\,r_k}.$$ Lemma~\ref{lem:depth} with $r=r_k$ gives a window $\{k-j',\dots,k-1\}$, $1\le j'\le k$, with more than $t_{j'}(r_k)$ misaligned iterations, so $\mcE_k$ occurs, contradicting $k\ge K_2'$.

\emph{(i).} Suppose $\varrho\le\frac12$ and that, on a path outside the null sets, $\|\nabla f(\BFX_k^0)\|\ge\epsilon$ for all $k\ge\bar k$, for some $\epsilon>0$ and $\bar k$. Decreasing $\epsilon$ and increasing $\bar k$ we may assume $\gamma_{inc}b_\epsilon\le\Delta_{\max}$ and $\bar k\ge\Ku$. Then~\eqref{eq:fixedpRadius} gives $$\Delta_k\ge \min\{\gamma_{dec}b_\epsilon,\Delta_{\bar k}\}\gamma_{dec}k^{-\varrho}$$ for all $k\ge K_2\vee(\bar k+1)$, so $$\sum_k\Delta_k^2\ge \left(\min\{\gamma_{dec}b_\epsilon,\Delta_{\bar k}\}\gamma_{dec}\right)^2\sum_{k\ge K_2\vee(\bar k+1)}k^{-2\varrho}=\infty,$$ contradicting Lemma~\ref{lem:squareSummable}.

\emph{(ii).} Let $\varrho<\frac12$; by (i), $\Teps<\infty$ for every $\epsilon>0$. Fix $\omega$. If $\Teps$ stays bounded as $\epsilon\downarrow0$ there is nothing to prove. Otherwise take $\epsilon\le1$ small enough that $\gamma_{inc}b_\epsilon\le\Delta_{\max}$, $\gamma_{dec}b_\epsilon\le\Delta_{\Ku}$ and $\Teps>K_2$. For $K_2\le k<\Teps$, \eqref{eq:fixedpRadius} with $\bar k=\Ku$ and $b_\epsilon=\ka\epsilon/\gamma_{dec}$ gives $\Delta_k\ge\gamma_{dec}\ka\epsilon\,k^{-\varrho}$, hence
\[
\SDelta\ge\sum_{k=K_2}^{\Teps-1}\Delta_k^2\ge\gamma_{dec}^2\ka^2\epsilon^2\int_{K_2}^{\Teps}x^{-2\varrho}\,dx=\gamma_{dec}^2\ka^2\epsilon^2\,\frac{\Teps^{1-2\varrho}-K_2^{1-2\varrho}}{1-2\varrho},
\]
so that $$\Teps^{1-2\varrho}\le\bigl((1-2\varrho)\gamma_{dec}^{-2}\ka^{-2}\SDelta+K_2^{1-2\varrho}\bigr)\epsilon^{-2},$$ which is the claim with $$\bar{R}:=\bigl((1-2\varrho)\gamma_{dec}^{-2}\ka^{-2}\SDelta+K_2^{1-2\varrho}\bigr)^{1/(1-2\varrho)}.$$ For the sample complexity, let $K_3:=K_2\vee\Kw$ and $Q_w:=\sum_{k<K_3}(2p+2)N_k$. For $K_3\le k<\Teps$, Lemma~\ref{lem:Nk} gives $$N_k\le n_W\lambda_k\Delta_k^{-2}\le n_W\lambda_{\Teps}\gamma_{dec}^{-2}\ka^{-2}\epsilon^{-2}k^{2\varrho},$$ so
\[
\Weps\le Q_w+(2p+2)\,n_W\lambda_{\Teps}\gamma_{dec}^{-2}\ka^{-2}\,\epsilon^{-2}\,\Teps^{1+2\varrho},
\]
and $$\epsilon^{-2}\Teps^{1+2\varrho}\le \bar{R}^{\,1+2\varrho}\epsilon^{-2-2(1+2\varrho)/(1-2\varrho)}=\bar{R}^{\,1+2\varrho}\epsilon^{-4/(1-2\varrho)},$$ while $\lambda_{\Teps}$ is polylogarithmic in $1/\epsilon$.
\end{proof}

\begin{remark}[the exponent]\label{rem:readingExponent}
For the ensembles of Section~\ref{sec:sketch}, $\dQ=e^{-c_{\sfQ}p}$, so $\log(1/\rho_{\sfQ})=c_{\sfQ}p-\log(e(1+c_\gamma)/c_\gamma)$ and
\[
\varrho=\frac{(1+c_\gamma)(1+\delta')\log(1/\gamma_{dec})}{c_{\sfQ}p-\log(e(1+c_\gamma)/c_\gamma)}=\mcO\Bigl(\frac1p\Bigr),\qquad \frac{2}{1-2\varrho}=2+\mcO\Bigl(\frac1p\Bigr).
\]
There are three consequences. First, a fixed subspace dimension is admissible on the pathwise route once $p$ exceeds the explicit threshold making $\varrho<\frac12$; the argument of this section then gives exponent $2+\mcO(1/p)$. 
Second, 
$p_k$ of order $\log k$ is the growth at which $k^{-\varrho_k}$ becomes a constant and the exponent is exactly $2$, which is Theorem~\ref{thm:iter}. 
Third, the threshold $\dQ<c_\gamma/(e(1+c_\gamma))$ comes from the crude Chernoff bound of Lemma~\ref{lem:chernoff}; replacing $(ej'\dQ/t)^t$ by the binomial Chernoff bound $\exp(-j'\,\mathrm{KL}(t/j'\|\dQ))$ extends the argument to $\dQ<c_\gamma/(1+c_\gamma)$, which for $\gamma_{inc}=1/\gamma_{dec}$ is $\dQ<\frac12$, the same threshold as Part~I.
\end{remark}

\begin{remark}[why expectation gets $\epsilon^{-2}$ with fixed $p$ and the path does not]\label{rem:whyLog}
The proof shows where the exponent comes from. Given the past, the depth of an excursion of the radius below $b_\epsilon$ has a geometric tail: $\mbP(\text{depth}\ge r)\le c\rho_{\sfQ}^{r/(1+c_\gamma)}$. In expectation this is harmless: $\mbE[\Delta_k^{-2}\ind\{k<\Teps\}]$ is bounded by a geometric series in $r$ as soon as $\gamma_{dec}^{-2}\rho_{\sfQ}^{1/(1+c_\gamma)}<1$, the mean depth is $\mcO(1)$, and one recovers $\mbE[\Teps]=\mcO(\epsilon^{-2})$; this is the mechanism behind the expected bounds of~\cite{blanchet2019convergence,chen2018stochastic} and of Part~I, whose Assumption~2.2 is the corresponding contraction condition for its inverse-radius sums. Along a path, however, the statement must hold for all $k$ at once, so what matters is not the mean depth but the maximum depth over the first $k$ iterations, which for a geometric tail is $\asymp\log k/\log(1/\rho_{\sfQ})$. That logarithm, exponentiated by $\gamma_{dec}$, is the factor $k^{-\varrho}$. Nothing is lost when the depth's logarithm is beaten by the growth of $p_k$, which is what $\sum_k(1-\bQk)<\infty$ achieves: then there are eventually no misaligned iterations at all, and the depth is zero.
\end{remark}


\section{Dimension dependence}\label{sec:dimension}

All statements about $d$, $p$ and the sketch scale are confined to this section; Sections~\ref{sec:accuracy}--\ref{sec:fixedp} are written in the named constants $\kQ,\aQ,\kedF,\kH,\kg,\kmf,\ksuc,\ka,\rF$. We take $p$ fixed for the accounting; with $p_k$ as in~\eqref{eq:pk} every power of $p$ below is a power of $\log(1/\epsilon)$, and $\aQ$ is evaluated at $p_0$ as in Proposition~\ref{prop:Haar}. We hold $(\Delta_0,\Delta_{\max},\mu,\kappa_s,\kH,\kedF)$ fixed in subspace coordinates, as discussed after Proposition~\ref{prop:reparam}, and drop absolute constants and the fixed problem constants $\kLg,\sigma_D,b_D,\sigma_0,\kappa_s,\Delta_{\max},\eta,\kappa_{\mathrm{fcd}},\gamma_{inc},\gamma_{dec}$ from the $\asymp$ relations. The random quantities $\SDelta$ and the burn-in indices also depend on $d$, $p$ and the scale; only the deterministic prefactors are tracked here, and the table takes $p_0=p$, since $\aQ$ is evaluated at $p_0$.

\paragraph*{The chain of constants.} From~\eqref{eq:constants}, \eqref{eq:ksuc} and~\eqref{eq:ka}, using $\ksuc^{-1}(p)\asymp\kmf(p)\vee\kH\vee\mu^{-1}$,
\begin{equation}
\kg(p)\asymp\sqrt p\,(\kQ^2+\kedF),\qquad \ka^{-2}(p)\asymp\aQ^{-2}\bigl(\sqrt p\,(\kQ^2+\kedF)+1\bigr)^2 ,\label{eq:kaDim}
\end{equation}
the additive $1$ coming from $\kH\ge1$ and $\mu^{-1}$, which do not rescale with the sketch. For the ensembles of Section~\ref{sec:sketch}, $\aQ^{-2}\asymp d/(p\kQ^2)$. On the sampling side, \eqref{eq:rF} gives
\begin{equation}
\frac1{\rF}\asymp\max\Bigl\{1,\ \Bigl(\frac{\kQ}{\kedF}\Bigr)^2,\ \frac{\kQ}{\kedF}\Bigr\},\label{eq:rFDim}
\end{equation}
and Theorem~\ref{thm:uniform} consumes $\lambda_k\ge2p/\rF$. Whether this is met by the burn-in index $k_\lambda$ or by choosing $\lambda_0\ge2p/\rF$, the sample size that the accuracy guarantee requires per iteration is of order $(p/\rF)\Delta_k^{-2}$, and this is the quantity to track: hidden in the index it reappears as a burn-in of order $\exp((p/\rF)^{1/(1+\epsilon_\lambda)})$, which for the unscaled sketch below is exponential in $d$. With $\Delta_k\ge\ka\epsilon$ before $\Teps$,
\begin{equation}
\Teps\lesssim\ka^{-2}\epsilon^{-2},\qquad N_k\lesssim\frac p{\rF}\ka^{-2}\epsilon^{-2},\qquad \Weps\asymp pN_k\Teps\lesssim\frac{p^2}{\rF}\ka^{-4}\epsilon^{-4},\label{eq:master}
\end{equation}
all up to polylogarithmic factors.

\paragraph*{Where the factor $p$ in $N_k$ comes from.} Lemma~\ref{lem:chained} controls $\sup_{\|\BFshat\|\le\Delta_k}|\Dbar_k(\BFshat,n)|$ at the scale $\kQ\Delta_k\sqrt{(p+\vartheta)/n}$, the $\sqrt p$ being the Gaussian width of the $p$-dimensional ball, and this must be brought below $\kedF\Delta_k^2$. The factor is intrinsic to a uniform certificate: the supremum of a linear functional over a $p$-ball is the norm of a $p$-vector of averaged noise, of order $\sqrt{p/d}$. It is what certifying accuracy at a point selected by the replications costs \emph{in the analysis}; a pointwise bound at an independently evaluated trial point does not pay it, and whether the algorithm pays it is a different question, since the factor enters a sufficient condition and no lower bound is given. In a subspace with $p_k=\mcO(\log k)$ the factor is polylogarithmic, so the certificate for reuse is essentially free there. In the full space, $p=d$, it is a factor $d$ in the sample size of our bound, or an exponential-in-$d$ burn-in if hidden in $k_\lambda$, so our bound for the reuse variant is one power of $d$ worse than for a fresh-stream variant while the saving in oracle calls per iteration is a factor $(2d+3)/(2d+2)$. We read this as a limitation of the certificate rather than as advice against reuse, which is what ASTRO-DF has always done in the full space. 

\paragraph*{Choosing the scale.} 
With $\kedF$ fixed and $\bar c:=\kQ^2$, \eqref{eq:kaDim} implies $\ka^{-2}\asymp d\,(\bar c+\kedF+p^{-1/2})^2/{\bar c}$, minimized at $\bar c=\kedF+p^{-1/2}$ with value $\asymp d\,(\kedF+p^{-1/2})$; lowering $\kedF$ therefore helps down to $\kedF\asymp p^{-1/2}$, beyond which the additive term dominates. 
The optimum is
\begin{equation}
\kedF\asymp p^{-1/2},\qquad\kQ\asymp p^{-1/4},\qquad\ka^{-2}\asymp d\,p^{-1/2},\qquad\frac1{\rF}\asymp p^{1/2}.\label{eq:balanced}
\end{equation}
The tolerance $\kedF$ is a parameter of the analysis, not of the algorithm: every result holds for each $\kedF\in(0,\ta)$, and a smaller value is paid for through~\eqref{eq:rFDim}, by the factor $1/\rF\asymp p^{1/2}$ in the sample size or by a later burn-in $k_\lambda$. 
In the deterministic analysis of~\cite{cartis2026note} the corresponding constant is $\kappa_d\asymp\sqrt p\kQ^2$ with no estimation term, and the same balance against the additive $1$ gives $\kQ\asymp p^{-1/4}$ there too; \cite[Remark~2.5]{cartis2026note} treats $p=\mcO(1)$ and so records only that scales below $d^{-1/2}$ gain nothing, which is correct when $p$ is not tracked.

\vspace{1em}
\begin{table}[h]\centering\small
\caption{Dimension dependence up to polylog factors, with $(\Delta_0,\Delta_{\max},\mu,\kappa_s,\kH,\kedF)$ fixed in subspace coordinates. Deterministic columns from~\cite{cartis2026note,chaudhry2026powell} with $\mcO(p)$ evaluations per iteration; stochastic columns from~\eqref{eq:master}. The unscaled sketch $\sqrt{d/p}\,\sfU_p$  has the same constants as an unscaled Gaussian or hashing transform, and $\sfS^\intercal/\sqrt d$ the same as $\sfU_p/\sqrt p$.}\label{tab:dim}
\begin{tabular}{llcccccc}
\toprule
& & & & \multicolumn{2}{c}{deterministic} & \multicolumn{2}{c}{stochastic, reuse}\\
\cmidrule(lr){5-6}\cmidrule(lr){7-8}
sketch $\sfQ_k$ & $\kedF$ & $\kQ$ & $\aQ$ & $\Teps\epsilon^2$ & evals$\cdot\epsilon^2$ & $\Teps\epsilon^2$ & $\Weps\epsilon^4$\\
\midrule
full space, $\sfI_d$ & $1$ & $1$ & $1$ & $d$ & $d^2$ & $d$ & $d^4$\\
unscaled, $\sqrt{d/p}\,\sfU_p$ & $1$ & $\sqrt{d/p}$ & $1$ & $d^2/p$ & $d^2$ & $d^2/p$ & $d^5/p$\\
orthonormal columns, $\sfU_p$ & $1$ & $1$ & $\sqrt{p/d}$ & $d$ & $pd$ & $d$ & $p^2d^2$\\
rescaled JLT, $\sfS^\intercal/\sqrt d$ & $1$ & $1/\sqrt p$ & $1/\sqrt d$ & $d$ & $pd$ & $pd$ & $p^4d^2$\\
rescaled JLT, $\sfS^\intercal/\sqrt d$ & $1/p$ & $1/\sqrt p$ & $1/\sqrt d$ & $d$ & $pd$ & $d$ & $p^3d^2$\\
balanced, $p^{-1/4}\sfU_p$ & $1/\sqrt p$ & $p^{-1/4}$ & $p^{1/4}/\sqrt d$ & $d/\sqrt p$ & $\sqrt p\,d$ & $d/\sqrt p$ & $p^{3/2}d^2$\\
\bottomrule
\end{tabular}
\end{table}

\paragraph*{Insights from the results} 
\begin{enumerate}
    \item In $d$, every subspace scale with $\kQ=\mcO(1)$ gives $\Teps=\Otilde(d\epsilon^{-2})$ and $\Weps=\Otilde(d^2\epsilon^{-4})$, against $d$ and $d^4$ in the full space: a factor $d^2$ in the sample complexity, one $d$ from the point count and one from the entropy of the trust region. 
    \item  The unscaled sketch $\sqrt{d/p}\,\sfU_p$  is worse than the full space by $d/p$ in iterations and in samples, and its sample size carries a factor $d$ through $1/\rF\asymp d/p$; this is the case in which hiding $2p/\rF$ in the burn-in index conceals an exponential-in-$d$ burn-in. The one-line change to orthonormal columns removes all of this.
    \item The Cartis--Roberts scale $\sfS^\intercal/\sqrt d$ shrinks the bias below the estimation error and, with $\kedF$ fixed, loses a factor $p$ in iterations relative to orthonormal columns; matching $\kedF\asymp\kQ^2$ restores $d$ at a cost in $1/\rF$. This is the stochastic form of Proposition~\ref{prop:reparam}: the scale is a parameter choice, and in the stochastic setting it must be made jointly with the accuracy demanded of the estimates. 
    \item The balanced scale~\eqref{eq:balanced} gains $\sqrt p$ in iterations over every other row, deterministically and stochastically. 
    \item The empirical preference for small $p$ observed in experiments is not explained by any row; every bound prefers large $p$ in $\Teps$ and penalizes it only through the polynomial prefactors of $\Weps$.
\end{enumerate}

\section{Concluding remarks}\label{sec:conclusion}
The paper set out to make the pathwise reading of a stochastic trust-region algorithm available, with its assumptions, outputs and limits stated plainly, to readers accustomed to the expectation reading. In that reading the requirements are summability of failure probabilities, obtained here by inflating the sample size logarithmically and growing the subspace dimension logarithmically; the device is the first Borel--Cantelli lemma, which needs no independence and is why common random numbers cause no trouble; the outputs are almost sure statements with random constants and random thresholds; and the limit is that nothing is said in expectation, though Section~\ref{sec:tails} shows the burn-in has an explicit, superpolynomially small tail, so finite-time high-probability statements come from the same computation. Within that reading we have shown that the sample size being a stopping time is a property to be used, by optional stopping, rather than worked around by a union bound, which is also where the pathwise analysis of~\cite{ha2025complexity} is sharpened here; that a data-dependent trial point calls for a uniform bound, not a pointwise one, while a trial point evaluated on its own stream does not; that at a fixed subspace dimension the eventual-accuracy argument gives exponent $2+\mcO(1/p)$, with the loss traced to the running maximum of a geometric depth; that the quadratic term of the model is a bounded perturbation in the first-order theory; and that the dimension dependence of the constants, once the sketch scale is understood as a parametrization, favors orthonormal columns and a balanced scale, with reuse costing, in our bounds, a factor $p$ that is polylogarithmic in a subspace and a genuine $d$ in the full space. Two things remain open: an explanation for the empirical preference for very small $p$, and whether the framework of Part~I can be extended to slowly varying constants so as to give the pathwise route expected bounds as well.

\subsection*{Acknowledgments}
This material was based upon work supported by National Science Foundation Grant CMMI-2226347, OAC-2410949 and Office of Naval Research Grant N000142412398. The first author gratefully acknowledges the discussions with Javier Madariaga regarding the finite-time high-probability bounds.

\bibliographystyle{plain}
\bibliography{refs}
\end{document}